\documentclass[a4paper]{amsart}
\usepackage{amsfonts}
\usepackage{amsmath}
\usepackage{mathrsfs}
\usepackage{graphicx}
\usepackage{amsmath,amsthm,amssymb,amscd}
\usepackage{color}
\usepackage{enumerate}

\usepackage{epsfig}

\theoremstyle{plain}

\newtheorem{Thm}{Theorem}[section]
\newtheorem{Lem}[Thm]{Lemma}
\newtheorem{Prop}[Thm]{Proposition}
\newtheorem{Cor}[Thm]{Corollary}
\theoremstyle{remark}
\newtheorem{Def}[Thm] {Definition}
\newtheorem{Rem}[Thm] {Remark}

\newtheorem{Ex}[Thm] {Example}

\newtheorem{thm}{Theorem}[section]

\newtheorem{defn}[thm] {Definition}

\theoremstyle{remark}

\theoremstyle{definition}

\newcommand{\eps}{\varepsilon}

\newcommand{\ra}{\rightarrow}

\newcommand{\htop}{h_{\text{top}}}

\long\def\begcom#1\endcom{}

\newcommand{\length}{\operatorname{\length}}

\def\length{\operatorname{length}}

\newcommand{\ul} {\underline}
\newcommand{\bl} {\begin{lemma}}
\newcommand{\el} {\end{lemma}}
\newcommand{\bt} {\begin{theorem}}
\newcommand{\et} {\end{theorem}}

\newcommand{\Tkk} {\mathcal T _{k+1}}

\newcommand{\limk}{\lim_{k \ra \infty}}

\newcommand{\Tk}{\mathcal T _k}

\newcommand{\bp}{\begin{proof}}
\newcommand{\ep}{\end{proof}}
\newcommand {\be}{\begin{equation}}
\newcommand  {\ee} {\end{equation}}
\newcommand  {\beq} {\begin{eqnarray*}}
\newcommand  {\eeq} {\end{eqnarray*}}

\newcommand  {\bd} {\begin{definition}}
\newcommand  {\ed} {\end{definition}}

\newtheorem{theorem}{Theorem}[section]
\newtheorem{lemma}[theorem]{Lemma}

\theoremstyle{definition}
\newtheorem{definition}[theorem]{Definition}

\theoremstyle{remark}

\numberwithin{equation}{section}

\begin{document}

\title[Lyapunov  `Non-typical' Points ] { Lyapunov `Non-typical' Points of  Matrix Cocycles 
 and Topological Entropy} 

\author[X. Tian] {Xueting Tian}
\address[X. Tian]{School of Mathematical Science,  Fudan University\\Shanghai 200433, People's Republic of China}
\email{xuetingtian@fudan.edu.cn}

\keywords{  Lyapunov Exponents; Cocycles; Topological Entropy;    Shadowing  Property;   Hyperbolic Systems}
\subjclass[2010] { 37H15; 37D20; 37D25; 37C45;  37C50; 
}
\maketitle

\def\abstractname{\textbf{Abstract}}

\begin{abstract}\addcontentsline{toc}{section}{\bf{English Abstract}}
It follows from Oseledec Multiplicative Ergodic Theorem (or Kingman's Sub-additional Ergodic Theorem) that the   set of `non-typical'  points for which the Oseledec averages of a given continuous cocycle  diverge has zero measure with respect to any   invariant probability measure.    In strong contrast, for any H$\ddot{o}$der  continuous cocycles over  hyperbolic systems, in this article we show that either all ergodic measures have same Maximal Lyapunov exponents or the set of Lyapunov `non-typical' points have full topological entropy and packing topological entropy. Moreover,  we give  an estimate of Bowen Hausdorff entropy from below.



\end{abstract}

\section{Introduction} \setlength{\parindent}{2em}

\subsection{Lyapunov Exponents}



 Let $f$ be an invertible map of a compact metric space $X$ and let $A:X\rightarrow GL (m,\mathbb{R})$  be a  continuous matrix function. One main object of interest is the asymptotic behavior
of the products of $A$ along the orbits of the
transformation $f$, called cocycle induced from $A$: for $n>0$
 $$A (x,n):=A (f^{n-1} (x))\cdots A (f (x))A (x), $$ and $$A (x,-n):=A (f^{-n} (x))^{-1}\cdots A (f^{-2} (x))^{-1} A (f^{-1} (x))^{-1}=A (f^{-n}x,n)^{-1}.$$

The Maximal Lyapunov exponent  (or simply, MLE) of $A:X\rightarrow GL (m,\mathbb{R})$ at one point $x\in X$ is defined as the limit
  $$\chi_{max} (A,x):=\lim_{n\rightarrow \infty}\frac1n{\log\|A (x,n)\|},
   $$ if it exists. In this case $x$ is called to be (forward) {\it Max-Lyapunov-regular} (simply, ML-regular). 
    Otherwise, $x$ is {\it Max-Lyapunov-irregular} (simply, ML-irregular,  or called Lyapunov `nontypical' point). By  Kingman's Sub-additional Ergodic Theorem (or Oseledec's Multiplicative Ergodic Theorem), for
    any invariant  measure $\mu$ and $\mu$ a.e. point $x$,  MLE always exists and the function $\chi_{max} (A,x)$  is $f$-invariant.
     Define $\chi_{max}(A,\mu)=\int \chi_{max} (A,x) d\mu.$
      In particular, for
    any ergodic  measure $\mu$ and $\mu$ a.e. point $x$,  MLE always exists and is constant, equal to $\chi_{max}(A,\mu)$. 
Let   $MLI (A,f)$ denote the set of all   ML-irregular points, called ML-irregular set.
 Then it is of zero measure for any  
  invariant measures.

In this paper we will mainly pay attention to the dynamical complexity of ML-irregular points, in the sense of entropy.

\subsection{Results}

 The set of all invariant
measures and the set of all ergodic invariant measures are denoted
by $\mathcal{M}_f(X)$ and $\mathcal{M}_f^e(X)$, respectively. An ergodic measure $\mu$ is called Lyapunov-minimizing for cocycle $A$, if $$\chi_{max}(\mu,A)=\inf_{\nu\in \mathcal{M}_f(X)} \chi_{max}(\nu,A).$$  An ergodic measure $\mu$ is called Lyapunov-minimizing for cocycle $A$,
   if $$\chi_{max}(\mu,A)=\inf_{\nu\in \mathcal{M}_f(X)} \chi_{max}(\nu,A).$$

 \begin{Thm}\label{Thm-LyapunovIrregular-HolderCocycle}
  Let $f$ be  a $C^{1}$   diffeomorphism of a compact Riemanian manifold $M$
 and $X\subseteq M$ be a topologically mixing locally maximal hyperbolic invariant subset.
     Let $A:X\rightarrow GL (m,\mathbb{R})$  be a H$\ddot{o}$der  continuous matrix function. Then either \\
     (1)  all ergodic measures supported on $X$ have same MLE w.r.t. $A$;  or \\
     (2)  ML-irregular set $MLI (A,f)\cap X$ has full topological entropy  of $X$ and  has full packing topological entropy of $X$ and moreover,  its Bowen Hausdorff entropy   can be estimated from below by $$\sup_{\mu \in \mathcal{M}^e_{f}(X)} \{ h_\mu(f) \,\,\,|\,\,\,\,\, \mu \text{ is not Lyapunov minimizing for } A\},$$ where $h_\mu(f)$ denotes the metric entropy of $\mu.$   In particular, if
 the  (unique) maximal entropy measure of $f|_X$ is not Lyapunov minimizing, then  $MLI (f)\cap X$ has  full Bowen Hausdorff entropy of $X$.

\end{Thm}

It is well-known that the maximal entropy measure of system restricted on any topologically mixing locally maximal hyperbolic set exists and unique, see \cite{Bowen2} (or \cite{DGS}).

\begin{Rem}  From \cite [Theorem 4]{Furman} we know    for any uniquely ergodic system, there exists some continuous matrix cocycle  whose Lyapunov-irregular set is a dense set of   second Baire category. So, it is admissible  to satisfy  that all ergodic measures have same Lyapunov spectrum and  simultaneously   Lyapunov-irregular points form a set of  second Baire category. On the other hand, if $m=1,$
the above phenomena of \cite{Furman} naturally does not happen for any dynamical system. Let us explain more precisely. If $m=1,$ the Lyapunov exponent can be written as Birkhoff ergodic average $$ \lim_{n\rightarrow +\infty}\frac1n\sum_{j=0}^{n-1}\phi(f^j(x))$$ where $\phi(x)=\log\|A (x)\|$ is a continuous function.  If all ergodic measures have same Lyapunov spectrum, then by Ergodic Decomposition theorem so do  all invariant measures and thus by weak$^*$ topology, the limit  $  \lim_{n\rightarrow +\infty}\frac1n\sum_{j=0}^{n-1}\phi(f^j(x))$  should exist at every point $x\in X$ and equal to the given spectrum. Moreover,   the case of $m=1$ is in fact to study Birkhoff ergodic average and it has been studied for systems with specification or its variants by many authors,   see \cite{Barreira-Schmeling2000,BarBook, Th2010,Tho2012} and reference therein.
As said in \cite{Bar-Gel} 
that the study of Lyapunov exponents lacks today a satisfactory general approach for non-conformal case, since a complete understanding is just known for some cases such as requiring a clear separation of Lyapunov directions or some number-theoretical properties etc, see \cite{Bar-Gel} and its references therein.  Our result gives one such  characterization for non-conformal maps on `non-typical' points of  Lyapunov exponents.

\end{Rem}

\begin{Rem}
 From \cite{Morris,Bousch-Jenkinson} and its references therein, we know  for generic continuous function,  there is a unique  maximizing (resp., minimizing, sometimes called Ergodic Optimization) measure with zero entropy.  Moreover, many people aim to show that  `most' functions are optimized by measures supported on a periodic orbit. For example,  from \cite{Yuan-Hunt} we know for hyperbolic systems (Axiom A), any Lipschitz function optimized by a periodic orbit measure can be perturbed to
be stably optimized by this periodic orbit measure;  and from \cite{Contr-L-T} for expanding systems, the subset of functions uniquely optimized by measures supported on a
periodic point is open and dense in the space of H$\ddot{o}$der continuous functions (in topology of H$\ddot{o}$der norm). These discussion are for the particular case of matrix functions with $m=1,$ and it is still an open question for the characterization of   Lyapunov maximizing (resp., minimizing) measure for matrix functions with $m\geq 2$. However, above analysis gives us some  positive direction to conjecture or believe that for `most' cocycles of $m\geq 2$, the minimizing measure should have metric entropy less than full entropy. If this is true, then for `most'  cocycles the Bowen Hausdorff entropy of Theorem \ref{Thm-LyapunovIrregular-HolderCocycle} is also full and same as another two entropy.  Moreover, on one hand,  it is unknown whether minimizing measure must exist naturally. If not, the Bowen Hausdorff entropy of Theorem \ref{Thm-LyapunovIrregular-HolderCocycle} is also full. On the other hand, we know that for hyperbolic system, there are infinite ergodic measure arbitrarily close to the full entropy (see entropy-dense property below), thus for a  matrix function  with $m\geq 2,$ if its minimizing measure exists but just forms a finite set, then the Bowen Hausdorff entropy of Theorem \ref{Thm-LyapunovIrregular-HolderCocycle} is  full. All in all, the better estimate of the Bowen Hausdorff entropy for Lyapunov `non-typical' points are based on the advance of Lyapunov Ergodic Optimization (there are few such results, just recently some people started to study, for example, see \cite{Bochi2014}).

\end{Rem}

As a particular case of Theorem \ref{Thm-LyapunovIrregular-HolderCocycle} we have a consequence on topological entropy for the derivative cocycle of
 all   hyperbolic systems.
 Let $MLI(f):=MLI(Df,f)$. It is called ML-irregular set of system $f$.

\begin{Thm}\label{Thm-horseshoe-LyapunovIrregular-Cocycle}
 Let $f$ be  a $C^{1+\alpha}$   diffeomorphism of a compact Riemanian manifold $M$
 and $X\subseteq M$ be a topologically mixing locally maximal invariant subset.
 Then either \\
     (1)  all ergodic measures  supported on $X$  have same MLE;  or \\
     (2)  ML-irregular set $MLI (f)\cap X$ has full topological entropy of $X$ and  has full packing topological entropy  of $X$ and moreover,   its Bowen Hausdorff entropy  can be estimated from below by
      $$\sup_{\mu \in \mathcal{M}^e_{f}(X)} \{ h_\mu(f) \,\,\,|\,\,\,\,\, \mu \text{ is not Lyapunov minimizing for } Df \}.$$
\end{Thm}

The motivations to obtain our main theorems are from some ideas of recent works of \cite{Feng-Huang,Kal,Tho2012,Dai}. In section \ref{section-prepare} we first give some preparation and then in section \ref{section-proof} we give the proof.

\section {Preliminaries}\label{section-prepare}

\subsection{ Entropy}

Let $(X,d)$ be a compact
metric space with Borel $\sigma-$algebra $\mathfrak{B}(X)$ and let
$f:X\rightarrow X$ be a continuous map. Let $\mathcal{M}(X)$ denote the space of all probability measures supported on $X$.  The set of all invariant
measures and the set of all ergodic invariant measures are denoted
by $\mathcal{M}_f(X)$ and $\mathcal{M}_f^e(X)$, respectively.
 For $x,y\in X$ and  $n\in \mathbb{N}$, let
$$d_n(x,y)=\max_{0\leq i\leq n-1}d(f^i(x),\,f^i(y)).$$   Let $x\in X$. The dynamical open ball $B_n (x,\varepsilon)$  and dynamical closed ball $\overline{B}_n (x,\varepsilon)$ are defined respectively  as
 $$B_n (x,\varepsilon):=\{y\in X|\, d_n(x,y) <\varepsilon\},\overline{B}_n (x,\varepsilon):=\{y\in X|\, d_n(x,y) \leq \varepsilon\}.$$  A set $S$ is $(n,\varepsilon)$-separated for $Z$ if $S\subset Z$ and $d_n(x,y)>\varepsilon$ for any $x,y\in S$ and $x\neq y$. A set $S\subset Z$ if $(n,\varepsilon)$-spanning for $Z$ if for any $x\in Z$, there exists $y\in S$ such that $d_n(x,y)\leq \varepsilon$.

 We have the following definition of entropy for compact set and thus the definition of entropy for a general subset.
\begin{defn}\label{compact-set}
For $E\subset X$ compact, we have the following Bowen's definition of topological entropy (see \cite{BowenEntropy}, c.f. \cite{Walter}).
\begin{equation}\label{Bowen-entropy}
  \htop(f,E)=\lim_{\eps\to 0}\lim_{n\to \infty}\frac{\log s_n(E,\eps)}{n}.
\end{equation}
For a general subset $Y\subset X$, define
\begin{equation}\label{entropy}
  \htop(f,Y)=\sup\{\htop(f,E):E\subset Y ~\textrm{is compact}\}.
\end{equation}
Finally, we put $\htop(f)=\htop(f,X)$. Since $X$ is a compact metric space, the definition depends only on the topology
on $X$, i.e. it is independent of the choice of metric defining the same topology on $X$.
\end{defn}

Let $\mu\in\mathcal{M}_f(X)$. Given $\xi=\{A_1,\cdots,A_k\}$ a finite measurable partition of $X$,
 i.e., a disjoint collection of elements of $\mathfrak{B}(X)$ whose
 union is $X$,  we define the entropy of $\xi$ by
 $$H_{\mu}(\xi)=-\sum_{i=1}^{k}\mu(A_i)\log \mu(A_i).$$
The metric entropy of $f$ with respect to $\xi$ is given by
 $$h_{\mu}(f,\,\xi)=\lim_{n\rightarrow\infty}\frac{1}{n}\log H_{\mu}(\bigvee_{i=0}^{n-1}f^{-i}\xi).$$
 The metric entropy of $f$ with respect to $\mu$ is given by
 $$h_{\mu}(f)=\sup_{\xi}h_{\mu}(f,\xi),$$
 where $\xi$ ranges over all finite measurable partitions of $X$.

\subsection{Bowen Hausdorff Entropy}

 Bowen also introduced another definition of entropy working for non-compact sets (see \cite{Bowen}). Some author call this notion \emph{Bowen Hausdorff entropy} (e.g. see \cite{Hass}).
\begin{defn}\label{non-compact-set}
For a general subset $E\subseteq X$, let $\mathcal {G}_{n}(E,\sigma)$ be the collection of all finite or countable covers of $E$ by sets of the form $B_{u}(x,\sigma)$ with $u\geq n$. We set
$$C(E;t,n,\sigma,f):=\inf_{\mathcal {C}\in \mathcal {G}_{n}(E,\sigma)}\sum_{B_{u}(x,\sigma)\in \mathcal {C}}e^{-tu}$$
and
$$C(E;t,\sigma,f):=\lim_{n\rightarrow\infty}C(E;t,n,\sigma,f).$$
Then
$$h_{top}^B(E;\sigma,f):=\inf\{t:C(E;t,\sigma,f)=0\}=\sup\{t:C(E;t,\sigma,f)=\infty\}$$
and the Bowen Hausdorff entropy of $E$ is
\begin{equation}\label{definition-of-topological-entropy}
  h_{top}^B(f,E):=\lim_{\sigma\rightarrow0} h_{top}(E;\sigma,f).
\end{equation}
Finally, we put $h_{top}^B(f)=h_{top}^B(f,X)$.
\end{defn}

It was proved by Bowen that $h_{top}^B(f)=\htop(f)$, however on subsets of $X$ we cannot guarantee this equality.
For example, if $A$ is a countable set then $h_{top}^B(f,A)=0$, while it may happen that $\htop(f,A)>0$.
In fact, it can be proved (cf. \cite{OZ}) that when invertible minimal dynamical system with finite entropy then for any $0\leq \alpha\leq \beta\leq \htop(f)$ there is plenty of compact sets $A$ with $h_{top}^B(f,A)=\alpha$
and $\htop(f,A)=\beta$. Then the cases $h_{top}^B(f,A)=\htop(f,A)$ are in some sense special.

{\rm Let $\mu\in \mathcal{M}(X)$.  The {\it measure-theoretical lower and
upper entropies} of $\mu$ are defined respectively by
$$\underline{h}_\mu(f)=\int \underline{h}_\mu(f,x) \;d\mu(x),\quad
\overline{h}_\mu(f)=\int \overline{h}_\mu(f,x) \;d\mu(x),$$ where
\begin{equation*}
\begin{split}
&\underline{h}_\mu(f,x)=\lim\limits_{\varepsilon\rightarrow
0}\liminf \limits_{n\rightarrow
+\infty}-\frac{1}{n}\log\mu(B_n(x,\varepsilon)),\\
&\overline{h}_\mu (f,x)=\lim\limits_{\varepsilon\rightarrow 0}\limsup
\limits_{n\rightarrow +\infty}-\frac{1}{n}\log\mu(B_n(x,\varepsilon)).
\end{split}
\end{equation*}

}

Brin and Katok \cite{BK} proved that for any $\mu\in M(X,f)$,
$\underline{h}_\mu(f,x)=\overline{h}_\mu(f,x)$ for $\mu$ a.e.
$x\in X$, and $\int \underline{h}_\mu(f,x)\; d\mu(x)=h_\mu(f)$.
So for  $\mu\in \mathcal{M}_f(X)$,
$$\underline{h}_\mu(f)=\overline{h}_\mu(f)=h_\mu(f).$$
From \cite[Proposition 1.2]{Feng-Huang} we know that if
$E\subseteq X$ is  non-empty and  compact, then
\begin{equation}\label{equation-bowenentropy}
h_{top}^B(f, E)=\sup \{\underline{h}_\mu(f): \mu \in \mathcal{M}(X), \; \mu(E)=1\}.
\end{equation}
 Remark that   one similar idea to estimate $h_{top}^B$ can be found in \cite{TV},  called entropy distribution principle.

\begin{Prop} \label{Prop2015-Pressure-distribution} (Entropy Distribution Principle)  
Let $f : X \mapsto X$ be a continuous transformation. Let $Z \subseteq X$ be an arbitrary Borel set. Suppose there exists $\varepsilon > 0$ and $s \geq 0$ such that one can find a sequence of Borel probability measures $\mu_k$, a constant $K> 0$, and a limit measure $\nu$ of the sequence $\mu_{k}$ satisfying $\nu(Z) > 0$ such that
\[
\limsup_{k \rightarrow \infty} \mu_{k} (B_n (x, \varepsilon)) \leq K \exp \{-ns \}
\]
for sufficiently large $n$ and every ball $B_n (x, \varepsilon)$ which has non-empty intersection with $Z$. Then $h_{top}^B(Z, \varepsilon,f) \geq s$.
\end{Prop}

\subsection{Packing topological entropy}
Let $E\subseteq X$. For $s\geq 0$, $N\in \mathbb{N}$ and $\varepsilon>0$, define
$$
P^s_{N,\varepsilon}(E)=\sup\sum_i
\exp(-sn_i),
$$
 where the supermum is taken over all finite or countable pairwise disjoint families $\{\overline{B}_{n_i}(x_i,\varepsilon)\}$ such that $x_i\in E$, $n_i\geq N$ for all $i$, where
 $$\overline{B}_n(x,\varepsilon):=\{y\in X:\; d_n(x,y)\leq \varepsilon\}.
 $$
 The quantity $P^s_{N, \varepsilon}(E)$ does not decrease as $N,\varepsilon$ decrease, hence the following
limits exist:
$$P^s_\varepsilon(E) = \lim_{N\to \infty} P^s_{N,\varepsilon}(E).$$
Define
$$
{\mathcal P}^s_\varepsilon(E)=\inf\left\{\sum_{i=1}^\infty P^s_\varepsilon(E_i):\; \bigcup_{i=1}^\infty E_i\supseteq E\right\}.
$$
Clearly, ${\mathcal P}^s_\varepsilon$ satisfies the following property:   if $E\subseteq \bigcup_{i=1}^\infty E_i$, then ${\mathcal P}^s_\varepsilon(E)\leq \sum_{i=1}^\infty {\mathcal P}^s_\varepsilon(E_i)$.

There exists a critical value
of the parameter s, which we will denote by $h_{top}^P(f,E,\varepsilon)$,
where ${\mathcal P}^s_\varepsilon(E)$ jumps from $\infty$ to $0$, i.e.
\[
{\mathcal P}^s_\varepsilon(E) = \left\{
\begin{array}{ll}
0, & s > h_{top}^P(f,E,\varepsilon),\\
\\
\infty,& s < h_{top}^P(f,E,\varepsilon).
\end{array}
\right.
\]
Note that  $h_{top}^P(f,E,\varepsilon)$ increases when $\varepsilon$ decreases.
 We  call
$$h_{top}^P(f, E):=\lim_{\varepsilon\to 0}h_{top}^P(f,E,\varepsilon)$$ the {\it packing topological entropy of $f$
restricted to $E$} or, simply, the {\it packing topological entropy
of $E$}, when there is no confusion about $f$. This quantity is defined in way which resembles the packing dimension. This definition is from \cite{Feng-Huang}.

From \cite[Proposition 2.1]{Feng-Huang} we know that
 for any $E\subseteq X$,
\begin{equation}\label{equation-entropydifferntrelation}
h_{top}^B(f, E)\leq h_{top}^P(f,E)\leq h_{top}(f,E);
\end{equation}
and  if $Z\subseteq X$ is  $f$-invariant and  compact, then \begin{equation}\label{equation-same-full-Packentropy}
h_{top}^B(f, Z)= h_{top}^P(f, Z)=h_{top}(f, Z).
\end{equation}
Moreover, from \cite[Proposition 1.3]{Feng-Huang} we know that if
$E\subseteq X$ is  non-empty and  compact, then
\begin{equation}\label{equation-Packentropy}
h_{top}^P(f, E)=\sup \{\overline{h}_\mu(f): \mu \in \mathcal{M}(X), \; \mu(E)=1\}.
\end{equation}


\subsection{Katok's definition  of metric entropy}

We use Katok's definition of metric entropy (see
\cite{Katok2}). It is defined for ergodic measure and equivalent to the classical one.
 Let $Z\subset X$. A set $S$ is $(n,\varepsilon)$-separated for $Z$ if $S\subset Z$ and $d_n(x,y)>\varepsilon$ for any $x,y\in S$ and $x\neq y$. A set $S\subset Z$ if $(n,\varepsilon)$-spanning for $Z$ if for any $x\in Z$, there exists $y\in S$ such that $d_n(x,y)\leq \varepsilon$.

Let $\mu\in \mathcal{M}^{e}_{f}(X)$.  For
$\varepsilon>0,\,\rho\in(0,1)$, let $N_l(\varepsilon,\,\rho)$ be the
minimal number of $\varepsilon-$balls $B_l(x,\,\varepsilon)$ in the
$d_l-$metric, which cover a set $Z\subseteq X$ of measure at least $1-\rho$.
 Define
 $$h_{\mu}(f,\varepsilon)= \liminf_{l\rightarrow\infty}\frac{\log N_l(\varepsilon,\,\rho)}{l},
\,\,h'_{\mu}(f,\varepsilon)= \limsup_{n\rightarrow\infty}\frac{\log
N_l(\varepsilon,\,\rho)}{l}.$$  Then
$$h_{\mu}(f)=\lim_{\varepsilon\rightarrow0}h_{\mu}(f,\varepsilon)
=h'_{\mu}(f,\varepsilon).$$

We need a following estimate for metric entropy.

\begin{Prop}\label{Prop2015-metric-entropy}
Let $\mu\in \mathcal{M}^{e}_{f}(X)$ and $\Gamma\subseteq X$ satisfy  $\mu(\Gamma)>0.$      For any $\delta>0,\varsigma>0$ $\varepsilon>0$, there is some positive integer $T_*$ (only dependent on $\delta$) such that for any  large $N>0$, there is $n\geq N$ and  $(n,\varepsilon)$-separated set $E_n\subseteq \Gamma$ such that \\
(1) for any $x\in E_n$, $f^n x\in \Gamma$;\\ 
(2) for any $x,y\in E_n$, $d(f^nx,y)<\delta$;\\ 
(3) $\frac 1n \log {T_*} \# E_n\geq  h_{\mu}(f,\varepsilon)  -(2+2\varsigma+h_{\mu}(f,\varepsilon))\varsigma. $ 

\end{Prop}

Before proving Proposition \ref{Prop2015-metric-entropy}, we need a following lemma.

\begin{Lem}\label{lem2015-recurrence}  Let $\mu\in \mathcal{M}^{e}_{f}(X)$ and  $\Gamma\subseteq X$ satisfy  $\mu(\Gamma)>0.$  Let $\delta>0$ and $\xi$ be a finite partition of $M$
with $\mbox{diam}\,\xi <\delta$ and $\xi>\{\Gamma,X\setminus
\Gamma\}$.   For $\rho>0,$ define
\begin{eqnarray*}
\Gamma_{s,\rho}&=&\big{\{}x\in \Gamma:
  \,\,f^m(x)\in \xi(x)~ \mbox{for some} ~
m\in[l, (1+{\rho})l],
 ~\mbox{for}~ l\geq s\big{\}}.
 \end{eqnarray*}
Then $\lim_{s\rightarrow +\infty} \mu (\Gamma_{s,\rho})=\mu (\Gamma). $
\end{Lem}

{\bf Proof.} This can be deduced from Birkhoff Ergodic Theorem. More precisely, let us first prove a following step.

\begin{Lem}\label{lem2015-recurrence-smallone}

 Let $\mu\in \mathcal{M}^{e}_{f}(X)$ and  $\Gamma\subseteq X$ satisfy  $\mu(\Delta)>0.$      For $\rho>0,$ define
\begin{eqnarray*}
\Delta_{s,\rho}&=&\big{\{}x\in \Delta:
  \,\,f^m(x)\in \Delta~ \mbox{for some} ~
m\in[l, (1+{\rho})l],
 ~\mbox{for}~ l\geq s\big{\}}.
 \end{eqnarray*}
Then $\lim_{s\rightarrow +\infty} \mu (\Delta_{s,\rho})=\mu (\Delta). $
\end{Lem}

 {\bf Proof.}  
   By Birkhoff Ergodic Theorem $\mu$ a.e. $y\in \Delta$,
$$ \lim_{n\rightarrow \infty} \frac 1n \sum_{i=0}^{n-1} \chi_{\Delta}(f^iy)=\mu(\Delta),$$
where $\chi_{\Delta}(y)$ denotes the  characteristic function on
 $\Delta.$  For given $\rho>0$, take $\varepsilon>0$ such that $0<\frac{1+\varepsilon}{1-\varepsilon}<1+\frac{3\rho}4$.  
 Then for $\mu$ a.e. $y\in \Delta$, there is $N_\varepsilon(y)>0$ such that for any $n\geq N_\varepsilon(y),$
\begin{eqnarray}\label{eqn2015-birkhoffestimate}
  1-\varepsilon< \frac 1{\mu(\Delta)} \frac 1n \sum_{i=0}^{n-1} \chi_{\Delta}(f^iy)<1+\varepsilon.
\end{eqnarray}
  Define $B_{s}=\{y\in \Delta:\,N_\varepsilon(y)\leq s\}$, then
$\lim_{s\rightarrow +\infty} \mu (B_{s})=\mu (\Delta). $
 Take $s\geq \frac 4 \rho,$ we claim $B_s\subseteq \Delta_{s,\rho}$. By contradiction,  there is $x\in B_s$ and  $l\geq s$
 such that $$f^m(x) \in X\setminus \Delta~ \mbox{for all} ~
m\in[l, (1+{\rho})l].$$ Then by the right estimate of (\ref{eqn2015-birkhoffestimate}) for $n=l$
$$  \frac 1{\mu(\Delta)} \frac 1{[(1+{\rho})l] } \sum_{i=0}^{[(1+{\rho})l]-1} \chi_{\Delta}(f^iy)= \frac 1{\mu(\Delta)} \frac 1{[(1+{\rho})l]} \sum_{i=0}^{l-1} \chi_{\Delta}(f^iy)<\frac 1{[(1+{\rho})l]} l(1+\varepsilon)$$
$$\leq \frac {l(1+\varepsilon)}{(1+{\rho})l-1}= \frac { (1+\varepsilon)}{(1+{\rho}) -\frac1l}\leq \frac { (1+\varepsilon)}{(1+{\rho}) -\frac\rho 4}<1-\varepsilon.$$  where $[a]$ denotes the integer part of $a$. This contradicts the left estimate in    (\ref{eqn2015-birkhoffestimate})  for $n=[(1+{\rho})l].$ \qed

\smallskip

Now we continue to prove Lemma \ref{lem2015-recurrence}. Let $\xi=\{C_1,C_2,\cdots, C_k\}\,(k\geq 1)$ and recall that $\xi>\{\Gamma,X\setminus
\Gamma\}$. Consider  $\Delta=C_i$ and by Lemma \ref{lem2015-recurrence-smallone} taking corresponding  $(C_i)_{s,\rho}$,
$$\lim_{s\rightarrow +\infty} \mu (\Gamma_{s,\rho})=\lim_{s\rightarrow +\infty} \mu (\cup_{i=1}^k \Gamma\cap (C_i)_{s,\rho})= \mu (\cup_{i=1}^k \Gamma\cap  C_i ) =\mu (\Delta). $$ \qed

{\bf Proof of Proposition \ref{Prop2015-metric-entropy}.} Take $\rho\in(0,1)$ such that  $1- \rho <\mu(\Gamma).$
Since $X$ is compact metric space, we can take  $\xi=\{\xi_1,\xi_2,\cdots,\xi_q\}$  to be a finite partition of $M$
with $\mbox{diam}\,\xi <\frac12\delta$ and $\xi>\{\Gamma,X\setminus
\Gamma\}$.  Fix a positive integer  $T_*\geq q=\# \xi.$
 Let $\Gamma_{s,\varsigma}$ be same as in Lemma \ref{lem2015-recurrence}. Then  $\mu(\Gamma_{s,\varsigma})\rightarrow \mu(\Gamma)\,(
as\,\,s\rightarrow+\infty)$ and thus we can take  sufficiently large $s$ such that
$\mu(\Gamma_{s,\varsigma})>  1-\rho.$ Let $D_l\subseteq \Gamma_{s,\varsigma}$ be an
$(l,\varepsilon)-$separated set of maximal cardinality; in other
words, the cover by $\varepsilon-$balls in the $d_l-$metric centered
at points in $D_l$ is a   cover of $\Gamma_{s,\varsigma}.$ Then we have
$$h_{\mu}(f,\varepsilon)\leq \liminf_{l\rightarrow+\infty}\frac1l\log \sharp
\, D_l.$$
So for given $N$ there exists $L>N$ such that,
$$h_{\mu}(f,\varepsilon)-\varsigma<\frac1L\log \sharp \,D_L, $$
and $$ {L\varsigma} <e^{ {L\varsigma} }.$$ Let $$G_l=\{x\in D_L: f^l(x)\in  \xi(x)\}, \,\,l=L,L+1,\cdots, [(1+\varsigma)L].$$
 By definition of $\Gamma_{L,\varsigma}$, from $ D_L\subseteq \Gamma_{L,\varsigma}$ we can take $n\geq N$ with $n\in[L,(1+\varsigma)L]$  such that
 \begin{eqnarray*}
 & & \# G_n\geq \frac1 {[(1+\varsigma)L]-L}\sum_{l=L}^ {[(1+\varsigma)L]-1}\# G_l \geq  \frac 1{[(1+\varsigma)L]-L}\# D_L \geq\frac 1{ \varsigma L }\# D_L \\
 & \geq &e^ {-L\varsigma+L(h_{\mu}(f,\varepsilon)-\varsigma)} =e^ {(-2\varsigma+h_{\mu}(f,\varepsilon))L }= e^ {(-2\varsigma+h_{\mu}(f,\varepsilon))n } e^ {(-2\varsigma+h_{\mu}(f,\varepsilon))(L-n) } \\
 &\geq&  e^ {(-2\varsigma+h_{\mu}(f,\varepsilon))n } e^ {-|(-2\varsigma+h_{\mu}(f,\varepsilon))(L-n)| }\geq e^ {(-2\varsigma+h_{\mu}(f,\varepsilon))n } e^ {-(2\varsigma+h_{\mu}(f,\varepsilon))\varsigma L }\\
  &\geq &  e^ {(-2\varsigma+h_{\mu}(f,\varepsilon))n } e^ {-(2\varsigma+h_{\mu}(f,\varepsilon))\varsigma n }=e^ {nh_{\mu}(f,\varepsilon)  -(2+2\varsigma+h_{\mu}(f,\varepsilon))\varsigma n }.
  \end{eqnarray*}  Let $F_{j}:=G_n\cap \xi_j,$ then we can take some $F_{j_0}$, denoted by $E_n,$ such that
  $$\# E_n=\max_{j=1}^q \# G_n\cap \xi_j\geq \frac 1 q \# G_n \geq \frac 1{T_*} e^ {nh_{\mu}(f,\varepsilon)  -(2+2\varsigma+h_{\mu}(f,\varepsilon))\varsigma n }.$$
 This gives the estimate (3).  (1) is from the definition of $G_n$ and $G_n\supseteq E_n$. By  the choice $\xi$ with $diam(\xi)<\frac12\delta$ and $E_n\subseteq \xi_{j_0}$, for any $x,y \in E_n$, we have
 $$d(f^nx,y)\leq d(f^nx,x)+d( x,y)\leq 2diam(\xi_{j_0})< \delta.$$  This is item  (2). Now we complete the proof. \qed

\medskip



\subsection{ (Exponential) Shadowing }



We begin with the classical shadowing property. Let $(X,d)$ be a compact metric space and $f:X \mapsto X$ be a homeomorphism.
A sequence $\{x_n\}_{n\in\mathbb{Z}}\subset M$ is called a $\delta$-pseudo-orbit of $f$ if
$$d(f(x_n),x_{n+1})<\delta~~\textrm{for any}~~n\in\mathbb{Z}.$$
Moreover, a $\delta$-pseudo-orbit $\{x_n\}_{n\in\mathbb{Z}}$ is $\tau$-shadowed by an orbit of $y\in M$ if
$$d(f^n(y),x_n)<\tau~~\textrm{for any}~~n\in\mathbb{Z}.$$
Finally, we say $f:M\to M$ satisfies shadowing property if for any $\tau>0$, there exists  $\delta>0$ such that any $\delta$-pseudo-orbit is $\tau$-shadowed by some orbit.


Given $x\in M$ and $n\in\mathbb{N}$, let
   $$\{x,\,n\}:=
 \{f^j (x)\,|\,\,j=0,\,1,\,\cdots,\,n-1\}.$$ In other words, $\{x,\,n\}$
represents the orbit segment from $x$ with length $n$. For a sequence of points
$\{x_i\}_{i=-\infty}^{+\infty}$ in $M$
 and a sequence of positive integers
 $\{n_i\}_{i=-\infty}^{+\infty}$, we call $\{x_i,\,n_i\}_{i=-\infty}^{+\infty}$
 a $\delta$-pseudo-orbit,
 if $$   d (f^{n_i} (x_i),\,x_{i+1})<\delta$$ for all $i$.

 Given $\lambda>0$ and $\tau>0,$ we call a point $x\,\in
M$  a  $ (\tau,\lambda)$-shadowing point for a pseudo-orbit
$\big{\{}x_i,\,n_i\big{\}}_{i=-\infty}^{+\infty},$  if
$$   d\big{ (}f^{c_i+j} (x),f^j (x_i)\big{)}<\tau\cdot e^{-\min\{j,n_i-j\}\lambda},$$ $\forall\,\,
j=0,\,1,\,2,\,\cdots,\,n_i-1$ and $\forall\,\, i \in \mathbb Z$, where
$c_i$
 is defined as
  \begin {equation} \label{eq:exponential-respective-time-squens}c_i=\begin{cases}
 0,&\text{for }i=0\\
 \sum_{j=0}^{i-1}n_j,&\text{for }i>0\\
 -\sum_{j=i}^{-1}n_j,&\text{for }i<0.
\end{cases}
 \end {equation}

%



Now we start to introduce  {\it  exponential  shadowing} property.
\begin{Def}\label{Def2015-Shadow} Let $\lambda>0$. 
$f$ is called to have  $\lambda$-exponential shadowing  property, if the
following holds: for any $\tau>0$ there exists  $\delta>0$ such that  for   any
$\delta$-pseudo-orbit $\{x_i,\,n_i\}_{i=-\infty}^{+\infty}$, 
there is a
$ (\tau,\lambda)$-shadowing point $x\in X $ for  $\{x_i,\,n_i\}_{i=-\infty}^{+\infty}$.\\
\end{Def}

Moreover, we say $f$ to have periodic $\lambda$-exponential shadowing  property, if it has  $\lambda$-exponential  shadowing  property and furthermore, if
$\{x_i,n_i\}_{i=-\infty}^{+\infty}$ in above definition is periodic, i.e., there
 exists an integer $m>0$ such that $x_{i+m}=x_i$ and $n_{i+m}=n_i$ for all i,
 then the shadowing point $x$ is periodic.

\begin{Rem}\label{Rem-exp-shadowing}

   If  a homeomorphism $f$ is topologically conjugated to a homeomorphism $g$ satisfying (periodic)  $\beta$-exponential shadowing property    for
 some $\beta>0$, and the inverse conjugation is $\gamma-$H$\ddot{\text{o}}$lder continuous, then it is not difficult to see that $f$   has  (periodic) $\beta\gamma$-exponential shadowing  property.
\end{Rem}

Remark that if above $\{x_i,n_i\}_{i=-\infty}^{+\infty}$ is periodic with $m=1,$ it is in fact the concept of {\it exponential closing}  introduced in \cite{Dai,Kal}.
 Exponential shadowing (or closing) plays important role in the estimate of Lyapunov exponents, for example, see \cite{Kal,WS,OT,STV}.

For convenience, we say the orbit segments $x,fx,\cdots,f^nx$ and $y,fy,\cdots,f^ny$ are exponentially $\tau$ close with exponent $\lambda$, meaning that $$d (f^i (x),f^i (y))<\tau e^{-\lambda\min\{i,n-i\}},0\leq i\leq n-1.$$

We will show that any hyperbolic set has exponential shadowing.

\begin{Prop}\label{Prop2015-LemShadow}
  Let $f:M\rightarrow M$ be a
$C^{1 }$ diffeomorphism, with an invariant hyperbolic locally maximal closed set $\Lambda.$ Then $f|_\Lambda$ has exponential shadowing property.
\end{Prop}

Before proving Proposition \ref{Prop2015-LemShadow},
 let us  recall Katok's shadowing lemma in the case of $C^1$ hyperbolic system
(c.f. \cite{Dai}, firstly proved in Chap. 5 \cite{Pollicott} for $C^{1+\alpha}$ non-uniformly hyperbolic case).
  Let $(\delta_{k})_{k=1}^{+
\infty}$ be a sequence of positive real numbers. Let
$(x_{n})_{n=-\infty}^{+ \infty}$ be a sequence of points in
$X$ for which there exists
a sequence $(s_{n})_{n=-\infty}^{+ \infty}$ of positive integers
satisfying:

\begin{enumerate}

\item[(i)] $\mid s_{n}-s_{n-1}\mid \leq 1, ~\forall n\in \mathbb{Z}$;

\item[(ii)] $d(fx_{n},x_{n+1})\leq \delta_{s_{n}}, ~\forall n\in \mathbb{Z}$;
\end{enumerate}
then we call $(x_{n})_{n=-\infty}^{+ \infty}$ a
$(\delta_{k})_{k=1}^{+ \infty}$  pseudo-orbit. Given $\tau>0$, a
point $x\in M$ is a
 $(\tau,\lambda)$-shadowing point for the $(\delta_{k})_{k=1}^{+ \infty}$
 pseudo-orbit if $d(f^{n}x,x_{n})\leq \tau \lambda_{s_{n}},~ \forall n\in
\mathbb{Z}$, where $\lambda_{k}=\lambda_{0}e^{-\lambda
k}$ and $\lambda_0$ is a constant only dependent on the system of $f$.

\begin{Lem}\label{Lem2015-Shadow}
(Katok's Shadowing lemma) Let $f:M\rightarrow M$ be a
$C^{1}$ diffeomorphism, with an invariant hyperbolic closed set $\Lambda.$ Then there exists $\theta_0>0$ and $\lambda_*>0$ such that for any   $  \tau \in(0,\theta_0),\lambda\in(0,\lambda_*)$,   there exists a
sequence $(\delta_{k})_{k=1}^{+ \infty}$ such that for any
$(\delta_{k})_{k=1}^{+ \infty}$ pseudo-orbit contained in $\Lambda$ there exists a unique
$(\tau,\lambda)$-shadowing
point in $M$.
\end{Lem}

{\bf Proof of Proposition \ref{Prop2015-LemShadow}.} Let $\theta_0>0$ and $\lambda_*>0$  same as in Lemma \ref{Lem2015-Shadow}. Fix  $\lambda\in(0,\lambda_*) $ and define corresponding sequence  $\lambda_{k}=\lambda_{0}e^{-\lambda
k}$ where $\lambda_0$ is a constant only dependent on the system of $f$.
 For $\tau>0$, 
by local maximal property of $\Lambda$, we can take  $\tau_1$ such that $0<\tau_1<\min\{ \tau, \theta_0\}$ and any  point $y$ with  $sup_{z\in Orb(y)}d(z,\Lambda)<\tau_1\lambda_0$ should be in $\Lambda$.   For $\tau_1$ and $\lambda$, using Lemma \ref{Lem2015-Shadow},  there exists a
sequence $(\delta_{k})_{k=1}^{+ \infty}$ such that for any
$(\delta_{k})_{k=1}^{+ \infty}$ pseudo-orbit contained in $\Lambda$ there exists a unique
$(\tau_1,\lambda)$-shadowing
point in $M$. By the choice of $\tau_1$, such showing point should be in $\Lambda.$

 Fix large $k$ such that  $\lambda_{k}=\lambda_{0}e^{-\lambda
k}<1.$  Take $\delta=\delta_k.$    For   any
$\delta$-pseudo-orbit $\{x_i,\,n_i\}_{i=-\infty}^{+\infty}$ contained in $\Lambda$,
 define $y_j=  f^{j-c_i} (x_i)$ and define $$ s_j=  {k+\min\{j-c_i,c_{i+1}-j\}},$$ for $c_i\leq j\leq c_{i+1}-1,\,\forall i\in \mathbb{Z}.$
  Then $\mid s_{j}-s_{j-1}\mid \leq 1, ~\forall j\in \mathbb{Z}$ and for $\,\forall i\in \mathbb{Z}$, $$d(f(y_j),y_{j+1})=0<\delta_{s_j},$$
 for $j$ with  $c_i  \leq j\leq c_{i+1}-2 $ and for $j=c_{i+1}-1$,
    $$d(f(y_j),y_{j+1})=d(f(y_{c_{i+1}-1}),y_{c_{i+1}})=d(f^{n_{i}}(x_{i}),x_{i+1}) <\delta=\delta_k=\delta_{s_j}.$$
This implies that   $\{x_i,\,n_i\}_{i=-\infty}^{+\infty}$ is a $(\delta_{k})_{k=1}^{+ \infty}$ pseudo-orbit contained in $\Lambda$. So there is
  a unique
$(\tau_1,\lambda)$-shadowing
point $y$ in $\Lambda$.  More precisely, for $c_i\leq j\leq c_{i+1}-1,\,  i\in \mathbb{Z},$
 $$ d(f^j(y),y_{j})<\tau_1 \lambda_{s_j}<\tau \lambda_0e^{-\lambda s_j}< \tau e^{-\lambda (s_j-k)}=\tau e^ {-\lambda\min\{j-c_i,c_{i+1}-j\} }.$$
 In other words, $$   d\big{ (}f^{c_i+j} (x),f^j (x_i)\big{)}<\tau\cdot e^{-\min\{j,n_i-j\}\lambda},$$ $\forall\,\,
j=0,\,1,\,2,\,\cdots,\,n_i-1$ and $\forall\,\, i \in \mathbb Z.$  \qed

\begin{Rem}\label{Rem-exp-shadowing-proof}
  One also can use another way to prove Proposition \ref{Prop2015-LemShadow}, since it is not difficult to see that
  $$\text{   Shadowing
+ Local product structure}
\Rightarrow \text{     Exponential  shadowing,}$$
 and it is known that every hyperbolic set has local product structure and every locally maximal hyperbolic set satisfies  shadowing property whose shadowing point is contained in the given set.
\end{Rem}

\begin{Ex}\label{exm-nonhyperbolic} {\it Non-hyperbolic systems with exponential shadowing: }
From \cite{Go} we know that  non-hyperbolic diffemorphism $f$ with $C^{1+Lip}$ smoothness, conjugated to a transitive Anosov diffeomorphism, exists even the conjugation and its inverse  is H$\ddot{\text{o}}$lder continuous. This example is non-hyperbolic but satisfies    exponential shadowing  property.



\end{Ex}





\subsection{Semi-Continuity of MLE 
 w.r.t. measures}

\begin{Prop}\label{Prop2015-Contiuity-LyaExp} Let $(X,d)$ be a compact metric space and $f:X \mapsto X$ be a homeomorphism.  Let $A:X\rightarrow GL (m,\mathbb{R})$  be a  continuous matrix function. Then
 $\chi_{max}(A,\mu)$ is upper semi-continuous with respect to $\mu \in \mathcal{M}_f(X).$

\end{Prop}

{\bf Proof.} Observe  that $$\chi_{max}(A,\mu)=\int  \chi_{max}(A,x)d\mu $$ $$
=\lim_{n\rightarrow +\infty}\int \frac1n{\log\|A (x,n)\|} d\mu=\inf_{n\geq 1}\int \frac1n{\log\|A (x,n)\|} d\mu.$$ The upper semi-continuity follows. \qed

\bigskip

Now  let's recall the entropy-dense property of Theorem 2.1 in    \cite{PS2005} (or see  \cite{PS}, also  see  \cite{EKW} for similar discussion).
    Roughly speaking, any invariant probability measure $\mu$ is the  limit of a sequence of ergodic measures
$\{\mu_n\}_{n=1}^{\infty}$ in weak$^*$ topology such that the entropy of $\mu$ is the limit of the entropies of  $\mu_n$.  We say $T$ has {\it entropy-dense} property, if for any $\nu\in \mathcal{M}_f(X)$, any neighborhood $G\subseteq \mathcal M (X)$ of $\nu$ and any $ h_* < h_\nu (f),$ there exists an ergodic
measure $\mu\in G\cap \mathcal{M}_f(X)$ such that   $h_\mu (f) > h_*.$

\begin{Prop}\label{Prop-Cor2015-Contiuity-LyaExp} Let $(X,d)$ be a compact metric space and $f:X \mapsto X$ be a homeomorphism.  Let $A:X\rightarrow GL (m,\mathbb{R})$  be a  continuous matrix function.  If there are two  ergodic measures with different MLE w.r.t. $A$, then \\
(1) $$h_{top}(f) =\sup_{\mu, \nu\in \mathcal{M}^e_{f}(X)} \{\max\{ h_\mu(f), h_\nu(f)\} |\,\,\chi_{max}(\mu,A)>  \chi_{max}(\nu,A)\};$$
 (2) If $f$ satisfies entropy-dense property, then  $$\sup_{\mu \in \mathcal{M}^e_{f}(X)} \{ h_\mu(f) |\,\,\chi_{max}(\mu,A)>\inf_{\nu\in \mathcal{M}_f(X)} \chi_{max}(\nu,A)\} $$
  $$=\sup_{\mu, \nu\in \mathcal{M}^e_{f}(X)} \{\min\{ h_\mu(f), h_\nu(f)\} |\,\,\chi_{max}(\mu,A)>  \chi_{max}(\nu,A)\}.$$ If further $\chi_{max}(A,\mu)$ is  lower semi-continuous with respect to $\mu \in \mathcal{M}_f(X),$ then they are equal to the full entropy, $h_{top}(f)$.

\end{Prop}

{\bf Proof.}
(1) Assume that  $$\sup_{\nu\in \mathcal{M}^e_f(X)} \chi_{max}(\nu,A) >\inf_{\nu\in \mathcal{M}^e_f(X)} \chi_{max}(\nu,A). $$  We only need to prove for any $\gamma>0,$ $$h_{top}(f) \leq \sup_{\mu, \nu\in \mathcal{M}^e_{f}(X)} \{\max\{ h_\mu(f), h_\nu(f)\} |\,\,\chi_{max}(\mu,A)>  \chi_{max}(\nu,A)\}+\gamma.$$
By classical variational principle (see Chapter 7 in \cite{Walter}),  there is an ergodic measure $\omega$ such that $h_\omega(f)\geq h_{top}(f)-\gamma.$
By assumption, take another $ \varpi\in \mathcal{M}^e_{f}(X)$ such that $\chi_{max}(\omega,A)\neq  \chi_{max}(\varpi,A)$.
Then $$h_{top}(f) \leq h_\omega(f) +\gamma \leq \max\{ h_\omega(f),h_\varpi(f)\}+\gamma $$
$$\leq \sup_{\mu, \nu\in \mathcal{M}^e_{f}(X)} \{\max\{ h_\mu(f), h_\nu(f)\} |\,\,\chi_{max}(\mu,A)>  \chi_{max}(\nu,A)\}+\gamma. $$

(2) The part of $\geq$ is trivial. So we only need to prove the part of $\leq.$  Suppose $\mu \in \mathcal{M}^e_{f}(X) $ satisfying
$\chi_{max}(\mu,A)>\inf_{\omega\in \mathcal{M}_f(X)} \chi_{max}(\omega,A) . $  Then there is $\omega\in \mathcal{M}_f(X)$ such that  $\chi_{max}(\mu,A)> \chi_{max}(\omega,A) . $  For $\gamma>0,$ take $\theta\in(0,1)$ close to 1 enough such that
$ \varpi=\theta \mu+(1-\theta)\omega$ has entropy larger than $h_\mu(f)-\gamma.$  Remark that $\chi_{max}(\varpi,A)< \chi_{max}(\mu,A)$ and take $\zeta<\chi_{max}(\mu,A)-\chi_{max}(\varpi,A) .$  By Proposition \ref{Prop2015-Contiuity-LyaExp} and the assumption of entropy-dense property, there is $\nu \in \mathcal{M}^e_{f}(X) $ close $\varpi$ enough such that $h_\nu(f)>h_\varpi(f)-\gamma>h_\mu(f)-2\gamma $ and
$$\chi_{max}(\nu,A)<\chi_{max}(\varpi,A)+\zeta < \chi_{max}(\mu,A). $$
Thus
 $$\sup_{\mu', \nu'\in \mathcal{M}^e_{f}(X)} \{\min\{ h_{\mu'}(f), h_{\nu'}(f)\} |\,\,\chi_{max}(\mu',A)>  \chi_{max}(\nu',A)\}$$$$\geq  \min\{ h_\mu(f), h_\nu(f)\} >h_\mu(f)-2\gamma.$$

  If further $\chi_{max}(A,\mu)$ is  lower semi-continuous with respect to $\mu \in \mathcal{M}_f(X),$ then we are going to prove   full entropy, $h_{top}(f)$. We only need to prove
 for any $\gamma>0,$ $$h_{top}(f)\leq \sup_{\mu \in \mathcal{M}^e_{f}(X)} \{ h_\mu(f) |\,\,\chi_{max}(\mu,A)>\inf_{\nu\in \mathcal{M}_f(X)} \chi_{max}(\nu,A)\}+2\gamma.$$
  By classical variational principle (see Chapter 7 in \cite{Walter}),  there is an ergodic measure ${\omega'}$ such that $h_{\omega'}(f)\geq h_{top}(f)-\gamma.$ If $\chi_{max}({\omega'},A)>\inf_{\nu\in \mathcal{M}_f(X)} \chi_{max}(\nu,A),$ then $$h_{top}(f)\leq h_{\omega'}(f)+\gamma\leq \sup_{\mu \in \mathcal{M}^e_{f}(X)} \{ h_\mu(f) |\,\,\chi_{max}(\mu,A)>\inf_{\nu\in \mathcal{M}_f(X)} \chi_{max}(\nu,A)\}+\gamma  .$$
Otherwise, $\chi_{max}({\omega'},A)=\inf_{\nu\in \mathcal{M}_f(X)} \chi_{max}(\nu,A).$ By assumption, take another $ {\varpi'}\in \mathcal{M}^e_{f}(X)$ such that $\chi_{max}({\omega'},A)\neq  \chi_{max}({\varpi'},A)$ and then $$\chi_{max}({\varpi'},A)>\inf_{\nu\in \mathcal{M}_f(X)} \chi_{max}(\nu,A)=\chi_{max}({\omega'},A).$$ By assumption of lower semi-continuity and entropy-dense,
 there is    ${\mu'} \in \mathcal{M}^e_{f}(X) $ close ${\varpi'}$ enough such that $h_{\mu'}(f)>h_{\varpi'}(f)-\gamma>h_{top}(f)-2\gamma $
  and
$$\chi_{max}(\mu',A)>\chi_{max}(\omega',A)=\inf_{\nu\in \mathcal{M}_f(X)} \chi_{max}(\nu,A). $$
Then $$h_{top}(f)\leq h_{\mu'}(f)+2\gamma\leq \sup_{\mu \in \mathcal{M}^e_{f}(X)} \{ h_\mu(f) |\,\,\chi_{max}(\mu,A)>\inf_{\nu\in \mathcal{M}_f(X)} \chi_{max}(\nu,A)\}+2\gamma  .$$
 \qed

\subsection{Lyapunov Exponents and Lyapunov Metric}

Suppose  $f:X\rightarrow X$ to be an invertible map on a compact metric space $X$ and $A:X\rightarrow GL (m,\mathbb{R})$ to  be a  continuous matrix function.

 \begin{Def}\label{def2015-Lyapunov-exp}
 For any $x\in X$ and any $0\neq v\in \mathbb{R}^m,$ define the Lyapunov exponent of vector $v$ at $x$,  $$\lambda (A,x,v):=\lim_{n\rightarrow +\infty}\frac1n{\log\|A (x,n)v\|},
   $$ if the limit exists. We say $x$ to be  (forward) {\it Lyapunov-regular} for $A$,  if $\lambda (A,x,v)$ exists for all vector
    $ v\in \mathbb{R}^m\setminus \{0\}.$ Otherwise, $x$ is called to be {\it Lyapunov-irregular} for $A$. Let $LI(A,f)$ denote the space of all Lyapunov-irregular points for $A$.
    \end{Def}

    By Oseledec's Multiplicative Ergodic theorem, for any invariant $\mu$ and $\mu$ a.e. $x$, the Lyapunov exponent
   $\lambda (A,x,v)$ exists at $x$ for all vectors $v\in \mathbb{R}^m.$



\bigskip

{\bf Oseledec Multiplicative Ergodic Theorem}  \cite [Theorem 3.4.4]{BP}: Let $f$ be an invertible ergodic measure-preserving transformation of a Lebesgue probability measure space $ (X,\mu).$ Let $A$ be a measurable cocycle whose generator satisfies $\log \|A^\pm (x)\|\in L^1 (X,\mu).$ Then there exist numbers $$\chi_1<\chi_2<\cdots<\chi_l,$$ an $f-$invariant set $\mathcal{R}^\mu$ with $\mu (\mathcal{R}^\mu)=1,$ and an $A-$invariant Lyapunov decomposition of $\mathbb{R}^m$ for $x\in \mathcal{R}^\mu,$
$$\mathbb{R}_x^m=E_{\chi_1} (x)\oplus E_{\chi_2} (x)\oplus \cdots E_{\chi_l} (x)$$ with $dim E_{\chi_i} (x)=m_i,$ such that
for any $i=1,\cdots,l$ and any $0\neq v\in E_{\chi_i} (x)$ one has
$$\lim_{n\rightarrow \pm\infty} \frac1n \log \|A (x,n)v\|=\chi_i$$
and $$\lim_{n\rightarrow \pm\infty} \frac1n \log det A (x,n) =\sum_{i=1}^lm_i\chi_i.$$

 \begin{Def}\label{Def-LyapunovExponents}
The numbers $\chi_1,\chi_2,\cdots,\chi_l$ are called the {\it Lyapunov exponents} of measure $\mu$ for cocycle $A$ and the dimension $m_i$ of the space $E_{\chi_i} (x)$ is called the {\it multiplicity} of the exponent $\chi_i.$ The collection of pairs $$Sp(\mu,A)=\{(\chi_i,m_i):1\leq i \leq l\}$$ is the {\it Lyapunov spectrum} of measure $\mu.$
\end{Def}

Remark that for any ergodic measure $\mu,$ all the points in the set $\mathcal{R}^\mu$ are  Lyapunov-regular.

We denote the standard scalar product in $\mathbb{R}^m$ by $<\cdot,\cdot>$. For a fixed $\epsilon>0$ and a regular point $x$ we introduce the {\it $\epsilon-$Lyapunov scalar product  (or metric)} $<\cdot,\cdot>_{x,\epsilon}$ in $\mathbb{R}^m$ as follows. For $u\in E_{\chi_i} (x),\, v\in E_{\chi_j} (x),\,i\neq j$ we define $<\cdot,\cdot>_{x,\epsilon}=0.$ For $i=1,\cdots,l$ and $u,v\in u\in E_{\chi_i} (x),$ we define
$$<\cdot,\cdot>_{x,\epsilon}=m\sum_{n\in\mathbb{Z}}<A (x,n)u,A (x,n)v>exp (-2\chi_in-\epsilon |n|).$$
Note that the series converges exponentially for any regular $x$. The constant $m$ in front of the conventional formula is introduced for more convenient comparison with the standard scalar product. Usually, $\epsilon$ will be fixed and we will denote $<\cdot,\cdot>_{x,\epsilon}$ simply by $<\cdot,\cdot>_{x}$
 and call it the {\it Lyapunov scalar product.} The norm generated by this scalar product is called the {\it Lyapunov norm}  and is denoted by $\|\cdot\|_{x,\epsilon}$ or $\|\cdot\|_{x}.$

 We summarize below some important properties of the Lyapunov scalar product and norm; for more details  see  \cite [\S 3.5.1-3.5.3]{BP}.  A direct calculation shows   \cite [Theorem 3.5.5]{BP} that for any regular $x$ and any $u\in E_{\chi_i} (x),\,\,\forall n\in \mathbb{Z},$
 \begin{eqnarray}\label{Lyapunov-norm}  exp (n\chi_i-\epsilon|n|)\|u\|_{x,\epsilon}\leq \|A (x,n)u\|_{f^nx,\epsilon}\leq exp  (n\chi_i+\epsilon|n|)\|u\|_{x,\epsilon},
 \end{eqnarray}
\begin{eqnarray}\label{Lyapunov-norm-Maximal-ightarrow}
 exp (n\chi-\epsilon|n|) \leq \|A (x,n)u\|_{f^nx\leftarrow x}\leq exp  (n\chi+\epsilon|n|),
 \end{eqnarray}
where $\chi=\chi_l$  is the maximal Lyapunov exponent and $\|\cdot\|_{f^nx\leftarrow x}$ is the operator norm with respect to the Lyapunov norms. It is defined for any matrix $A$ and any regular points $x,y$ as follows:
$$\|A\|_{y\leftarrow x}=\sup\{\|Au\|_{y,\epsilon}\cdot \|u\|^{-1}_{x,\epsilon}:0\neq u\in \mathbb{R}^m\}.$$

We emphasize that, for any given $\epsilon>0,$ Lyapunov scalar product and Lyapunov norm are defined only for regular points with respect to the given measure. They depend only measurably on the point even if the cocycle is H$\ddot{\text{o}}$lder. Therefore, comparison with the standard norm becomes important. The uniform lower bound follow easily from the definition: $$\|u\|_{x,\epsilon}\geq \|u\|.$$ The upper bound is not uniform, but it changes slowly along the regular orbits (\cite{BP}, Prop. 3.5.8): there exists a measurable function $K_\epsilon (x)$ defined on the set of regular points $\mathcal{R}^\mu$ such that
\begin{eqnarray}\label{eq-different-norm-estimate}
\|u\|\leq \|u\|_{x,\epsilon}\leq K_\epsilon (x)\|u\|\,\,\,\,\,\,\,\,\,\forall x\in\mathcal{R}^\mu,\,\,\forall u\in  \mathbb{R}^m
 \end{eqnarray}
\begin{eqnarray}\label{eq-estimate-K-epsilon}
  K_\epsilon (x)e^{-\epsilon n}\leq K_\epsilon (f^nx)  \leq K_\epsilon (x)e^{\epsilon n}\,\,\,\,\,\,\,\,\,\forall x\in\mathcal{R}^\mu,\,\,\forall n\in  \mathbb{Z}.
 \end{eqnarray}
 These estimates are obtained in \cite{BP} using the fact that $\|u\|_{x,\epsilon}$ is {\it tempered}, but they can also be checked directly using the definition of  $\|u\|_{x,\epsilon}$ on each Lyapunov space and noting that angles between the spaces change slowly.

 For any matrix $A$ and any regular points $x,y,$ inequalities  (\ref{eq-different-norm-estimate}) and  (\ref{eq-estimate-K-epsilon}) yield
\begin{eqnarray}\label{eq-estimate-norm-K-epsilon}
  K_\epsilon (x)^{-1}\|A\| \leq \|A\|_{y\leftarrow x}\leq K_\epsilon (y) \|A\|.
 \end{eqnarray}

When $\epsilon$ is fixed we will usually omit it and write $K (x)=K_\epsilon (x).$ For any $l>1$ we also define the following sets of regular points
\begin{eqnarray}\label{eq-estimate-measure-Pesinblock}
  \mathcal{R}^\mu_{\epsilon,l}=\{x\in \mathcal{R}^\mu: \,\,K_\epsilon (x)\leq l\}.
 \end{eqnarray}
Note that $\mu (\mathcal{R}^\mu_{\epsilon,l} )\rightarrow 1$ as $l\rightarrow \infty.$ Without loss of generality, we can assume that the set $\mathcal{R}^\mu_{\epsilon,l}$ is compact and that Lyapunov splitting and Lyapunov scalar product are continuous on $\mathcal{R}^\mu_{\epsilon,l}.$ Indeed, by Luzin's theorem we can always find a subset of $\mathcal{R}^\mu_{\epsilon,l}$ satisfying these properties with arbitrarily small loss of measure (for standard Pesin sets these properties are automatically satisfied).

\subsection{Estimate of the norm of H$\ddot{\text{o}}$lder cocycles}\label{Estimate-norm}

Let us recall some useful lemmas   as follows. Firstly let us recall a general estimate of the norm of $A$ along any orbit segment close to a regular one\cite{Kal}.

 \begin{Lem}\label{Lem-simple-estimate-Lyapunov}   \cite [Lemma 3.1] {Kal}
Let $A$ be an $\alpha-$H$\ddot{\text{o}}$lder cocycle ($\alpha>0$) over a continuous map $f$ of a compact metric space $X$ and let $\mu$ be an ergodic measure for $f$ with the largest Lypunov exponent $\chi.$ Then for any positive $\lambda$ and $\epsilon$ satisfying $\lambda>\epsilon/ \alpha$ there exists $c>0$ such that for any $n\in\mathbb{N}$, any regular point $x$ with both $x$ and $f^nx$ in $\mathcal{R}^\mu_{\epsilon,l}$, and any point $y\in X$ such that the orbit segments $x,fx,\cdots,f^nx$ and $y,fy,\cdots,f^n (y)$ are exponentially $\tau$ close with exponent $\lambda$ for some $\tau>0$ we have
\begin{eqnarray}\label{eq-estimate-simple-Lyapunov-1}
  \|A (y,n)\|_{f^nx\leftarrow x}\leq e^{cl\tau^\alpha}e^{n (\chi+\epsilon)}\leq e^{2n\epsilon+cl\tau^\alpha}\|A (x,n)\|_{f^nx\leftarrow x}
 \end{eqnarray}
 and
 \begin{eqnarray}\label{eq-estimate-simple-Lyapunov-2}
  \|A (y,n)\| \leq l^2 e^{cl\tau^\alpha}e^{n (\chi+\epsilon)}\leq l^2 e^{2n\epsilon+cl\tau^\alpha}\|A (x,n)\|.
 \end{eqnarray}
 The constant $c$ depends only on the cocycle $A$ and on the number $ (\alpha\lambda-\epsilon).$

\end{Lem}

\begin{Lem}\label{Lem-New-simple-estimate-Lyapunov-new}
Let $A$ be an $\alpha-$H$\ddot{\text{o}}$lder cocycle  ($\alpha>0$) over a continuous map $f$ of a compact metric space $X$ and let $\mu$ be an ergodic measure for $f$ with the largest Lypunov exponent $\chi.$ Then for any positive $\lambda$ and $\epsilon$ satisfying $\lambda>\epsilon/ \alpha$ there exists $\tau>0$ such that for any $n\in\mathbb{N}$, any regular point $x$ with both $x$ and $f^nx$ in $\mathcal{R}^\mu_{\epsilon,l}$, and any point $y\in X$ such that the orbit segments $x,fx,\cdots,f^nx$ and $y,fy,\cdots,f^n (y)$ are exponentially $\tau$ close with exponent $\lambda$ for some $\tau>0$ we have
 \begin{eqnarray}\label{eq-estimate-simple-Lyapunov-new}
  \|A (y,n)\| \leq l^2 e^{l}e^{n (\chi+\epsilon)}\leq l^2e^l e^{2n\epsilon} \|A (x,n)\|.
 \end{eqnarray}

\end{Lem}

{\bf Proof.} For Lemma \ref{Lem-simple-estimate-Lyapunov}, let $\tau>0$ small enough such that $$c\tau^\alpha<1.$$ Then the   estimate (\ref{eq-estimate-simple-Lyapunov-new}) is obvious from Lemma \ref{Lem-simple-estimate-Lyapunov}. \qed

\bigskip

Another lemma is to estimate the growth of vectors in a ceratin cone $K\subseteq \mathbb{R}^m$ invariant under $A (x,n)$ \cite{Kal}. Let $x$ be a point in $\mathcal{R}^\mu_{\epsilon,l} $ and $y\in X$ be a point such that the orbit segments $x,fx,\cdots,f^nx$ and $y,fy,\cdots,f^ny$ are exponentially $\tau$ close with exponent $\lambda.$ We denote
$x_i=f^ix$ and $y_i=f^iy,\,i=0,1,\cdots,n.$ For each $i$ we have orthogonal splitting $\mathbb{R}^m=E_i\oplus F_i,$ where $E_i$ is the Lyapunov space at $x_i$ corresponding to the largest Lyapunov exponent $\chi$ and $F_i$ is the direct sum of all other Lyapunov spaces at $x_i$ corresponding to the Lyapunov exponents less than $\chi.$ For any vector $u\in \mathbb{R}^m$ we denote by $u=u'+u^\perp$ the corresponding splitting with $u'\in E_i$ and $u^\perp\in F_i;$ the choice of $i$ will be clear from the context. To simplify notation, we write $\|\cdot\|_i$ for the Lyapunov norm at $x_i$. For each $i=0,1,\cdots,n$ we consider cones $$K_i=\{u\in \mathbb{R}^m:\,\|u^\bot\|_i\leq \|u'\|_i\}\,\,\,\,\text{ and } \,\,\,\,K_i^\eta=\{u\in \mathbb{R}^m:\,\|u^\bot \|_i\leq  (1-\eta)\|u'\|_i\}$$
with $\eta>0$. Remark that \begin{eqnarray}\label{eq-u-u'}
\|u\|_i\geq \|u'\|_i\geq\frac1{\sqrt2}\|u\|_i.
\end{eqnarray}
We will consider the case when $\chi$ is not the only Lyapunov exponent of $A$ with respcet to $\mu.$ Otherwise $F_i=\{0\}, K_i^\eta=K_i=\mathbb{R}^m$, and the argument becomes simpler. Recall that $\epsilon<\epsilon_0=\min\{\lambda\alpha, (\chi-\nu)/2\},$ where $\nu<\chi$ is the second largest Lyapunov exponent of $A$ with respect to $\mu.$

\begin{Lem}\label{Lem-simple-estimate-Lyapunov-2}  \cite [Lemma 3.3]{Kal}
In the notation above, for any regular set $\mathcal{R}^\mu_{\epsilon,l}$, there exist $\eta,\tau>0$ such that if  $x,f^nx\in\mathcal{R}^\mu_{\epsilon,l}$  and   the orbit segments $x,fx,\cdots,f^nx$ and $y,fy,\cdots,f^n (y)$ are exponentially $\tau$ close with exponent $\lambda$, then for every $i=0,1,\cdots,n-1$ we have $A (y_i) (K_i)\subseteq K_i^\eta$ and  $\| (A (y_i)u)'\|_{i+1}\geq e^{\chi-2\epsilon}\|u'\|_i$ for any $u\in K_i.$

\end{Lem}

\section {Proof}\label{section-proof}


\subsection{ML-irregular set and Topological Entropy}

\begin{Thm}\label{Thm-Entropy-Maximal-LyapunovIrregular-HolderCocycle}
Let $(X,d)$ be a compact metric space and $f:X \mapsto X$ be a topologically mixing  homeomorphism with   exponential shadowing
 property and let $A:M\rightarrow GL(m,\mathbb{R})$ be a H$\ddot{\text{o}}$lder continuous matrix function. Then
 either \\
     (1)  all ergodic measures   have same MLE  w.r.t. $A$;  or \\
     (2)  ML-irregular set $MLI (A,f) $ has full topological entropy    and  has full packing topological entropy   and moreover,  its Bowen Hausdorff entropy   can be estimated from below by $$\sup_{\mu \in \mathcal{M}^e_{f}(X)} \{ h_\mu(f) \,\,\,|\,\,\,\,\, \mu \text{ is not Lyapunov minimizing for } A\}.$$



\end{Thm}

{\bf Proof.}
 Let $A$ be   $\alpha-$H$\ddot{\text{o}}$lder continuous and let $C:=\max_{x\in X}\|A^\pm(x)\|$. Firstly, let us to show the estimate of entropy by metric entropy of two measures with different Lyapunov exponents.
Let $\mu_1$ and $\mu_2$ be two ergodic measures  such that  $$\int \chi_{max}(x)d\mu_1>\int \chi_{max}(x)d\mu_2.$$
   Let $h^*=\min\{h_{\mu_1}(f),\,h_{\mu_2}(f)\},\,H^*=\max\{h_{\mu_1}(f),\,h_{\mu_2}(f)\}.$

   \begin{Prop}\label{prop-entropy-twomeasure} $$ \min\{ h_{top}(MLI(A,f),f) , h_{top}^P(MLI(A,f),f)\}\geq
H^* ,\,\,\,h_{top}^B(MLI(A,f),f)\geq
h^* .$$
\end{Prop}

{\bf Proof.}
  From (\ref{equation-entropydifferntrelation}),  $h_{top}^P(MLI(A,f),f)\leq h_{top}(MLI(A,f),f)$,  thus by Proposition \ref{Prop-Cor2015-Contiuity-LyaExp} we only  need to show for any $\gamma>0,$ $$h_{top}^P(MLI(A,f),f)\geq
H^*-5\gamma,\,\,h_{top}^B(MLI(A,f),f)\geq
h^*-4\gamma.$$ Now fix $\gamma>0.$

\textbf{Step 1}. Choice of  separated sets.

  If let $a=\int \chi_{max}(x)d\mu_1$ and $b=\int \chi_{max}(x)d{{\mu_2}}$, then $a-b>0$.
 Let $\lambda $ be the positive number in the definition of exponential shadowing.
  Take  $\epsilon\in(0,\frac{a-b}8)$  satisfying $\lambda>\epsilon/ \alpha$.

  Fix  a number $\rho\in(0,1).$
    By Katok's definition of metric entropy,  take $\varepsilon>0$ such that
\begin{eqnarray}\label{eq-Katokentropy-choice} h_{\mu_i}(f,4\varepsilon)=\liminf_{n\rightarrow+\infty}\frac1n\log
N^{\mu_i}_n(4\varepsilon,\rho)
>h_{\mu_i}(f)-\gamma,\,\,\,i=1,2.
\end{eqnarray}
 Take  $\varsigma>0$ small enough such that \begin{eqnarray}\label{eq-Katokentropy-choice-2}
  h_{\mu_i}(f,4\varepsilon)  -(2+2\varsigma+h_{\mu_i}(f,4\varepsilon))\varsigma\geq h_{\mu_i}(f,4\varepsilon)  -\gamma,\,\,i=1,2.
  \end{eqnarray}

 For the measures ${\mu_1}$ and ${{\mu_2}},$ take $l$ large enough such that $${\mu_1}(\mathcal{R}^{\mu_1}_{\epsilon,l})>1-\rho,\,\,\,{\mu_2}(\mathcal{R}^{\mu_2}_{\epsilon,l})>1-\rho.$$
 Take $\eta>0,\tau\in(0,
 \varepsilon)$ small enough such that it is applicable to Lemma \ref{Lem-New-simple-estimate-Lyapunov-new}
    and Lemma \ref{Lem-simple-estimate-Lyapunov-2}.  For $\tau,\,\lambda,$ by exponential shadowing there is $ \delta>0$  such that  for   any
$\delta$-pseudo-orbit $\{x_i,\,n_i\}_{i=-\infty}^{+\infty}$,  
there is a
$ (\tau,\lambda)$-shadowing point $x\in X $ for  $\{x_i,\,n_i\}_{i=-\infty}^{+\infty}$.

   By compactness of $X$, we can choose a finite open cover $\{U_i\}_{i=1}^q$ with $diam(U_i)<\delta.$ Define $$T_{i,j}=\min\{n\,\,|\,\,U_i\cap f^{-n-p}U_j\neq \emptyset,\,\forall \,p\geq 0\}.$$ Since $f$ is topologically mixing, $T_{i,j}$ is well-defined. Take $N=\max\{T_{i,j}\,|\,\,1\leq i,j\leq q\}.$  Then for any $n\geq N$ and any $1\leq i,j\leq q$, \begin{eqnarray}\label{eq-mixing}
    U_i\cap f^{-n}U_j\neq \emptyset.
    \end{eqnarray}

Let $\mu=\mu_i$ and $\Gamma= \mathcal{R}^{\mu_i}_{\epsilon,l} $ in Proposition \ref{Prop2015-metric-entropy}.        For above $\delta>0,\varsigma>0$ $\varepsilon>0$,  there is $T_*>0$ such that for any large $N_*>0$, there is ${n_i}\geq N_*$ and  $({n_i},4\varepsilon)$-separated set $E^{i}_{n_i}\subseteq \mathcal{R}^{\mu_i}_{\epsilon,l}$ such that
for each  $i=1,2$,  \\
(1) for any $x\in E^i_{n_i}$, $f^{n_i} x\in \mathcal{R}^{\mu_i}_{\epsilon,l}$ and  for any $x,y\in E^i_{n_i}$, $d(y,f^{n_i}x)< \delta;$\\
(2) $\frac 1{n_i} \log T_*\# E^i_{n_i}\geq h_{\mu}(f,\varepsilon)  -(2+2\varsigma+h_{\mu}(f,\varepsilon))\varsigma. $\\ 
By (\ref{eq-Katokentropy-choice}) and (\ref{eq-Katokentropy-choice-2}), $\frac 1{n_i} \log T_*\#  E^i_{n_i}\geq h_{\mu_i}(f)-2\gamma,\,\,i=1,2. $ We emphisize  that
$n_1$ and $n_2$ can be chosen arbitririly large.  So we can take $n_1,n_2$ also  satisfy  following estimates:   $n_1 \gg N$ such that
 \begin{eqnarray}\label{eq-E-choice}n_1\gamma>(h^*-3\gamma)N
 ,\,\,\,\,2<  e ^{ n_1 \epsilon},\,\,\,\,\,\,\,\frac 1{n_1} \log \#  E^1_{n_1}\geq h_{\mu_1}(f)-3\gamma
 \end{eqnarray}
 and  $n_2\gg N$ 
  large enough such that
     \begin{eqnarray}\label{eq-E-choice-2} n_2\gamma>(h^*-3\gamma)N
     ,\,\,\,\, l^2 e^{l}<e^{ n_2\epsilon},\,\,\,\,\,\,\frac 1{n_2} \log \#  E^2_{n_2}\geq h_{\mu_2}(f)-3\gamma.\end{eqnarray}

Let $g:\mathbb{N}\rightarrow \{1,2\}$ be given by $g(k)=(k+1)(\text{mod }2)+1$ and let
$$\mathcal S_k:=E^{g(k)}_{n_{g(k)}},\,n_k=n_{g(k)},\,\mu_k=\mu_{g(k)}.$$
From above construction,   we have
\begin{eqnarray}\label{eq-Yk-choice}
\# \mathcal S_k\geq exp(h_{\mu_{k}}(f)-3\gamma)n_{k};
\end{eqnarray}
and \begin{eqnarray}\label{eq-Yk-choice-2}
\text{ for any } x\in  \mathcal S_k,  f^{n_{k}} x\in \mathcal{R}^{\mu_k}_{\epsilon,l} \text{  and  for any }  x,y\in  \mathcal S_k, d(y,f^{n_{k}}x)<\delta.
\end{eqnarray}

\textbf{Step 2}. {Construction of the fractal F.}

Let us choose a sequence with $N_0 = 0$ and $N_k$ increasing to $\infty$ sufficiently quickly so that
\begin{equation} \label{f.1}
\limk \frac{n_{k+1}}{N_k} = 0,  \limk \frac{N_1 n_1 + \ldots + N_k n_k+kN}{N_{k+1}} = 0.
\end{equation}

Let $\ul x_i =(x_1^i, \ldots, x_{N_i}^i) \in S_i^{N_i}$. For any $(\ul x_1, \ldots, \ul x_k) \in S_1^{N_1} \times \ldots \times S_k^{N_k}$, by the exponential shadowing property, we have $B(\ul x_1, \ldots \ul x_k) :=$
\[
 \bigcap_{i=1}^k \bigcap_{j=1}^{N_i} f^{-\sum_{l=0}^{i-1} N_l n_l-(i-1)N - (j-1) n_i}\bigcap_{t=0}^{n_i-1}f^{-t}B ( f^t x^i_j, \tau e^{\lambda \min\{t,n_i-t\}}) \neq \emptyset.
\]
Let us explain this  emptiness   more precisely. Fix $(\ul x_1, \ldots, \ul x_k) \in S_1^{N_1} \times \ldots \times S_k^{N_k}$.  Recall above  finite open cover $\{U_t\}_{t=1}^q$ and we can take $U_{t_i}$ and $U_{t_{i+1}}$ such that $f^{n_i}x_{N_i}^i\in U_{t_i}$ and $ x_{1}^{i+1}\in U_{t_{i+1}}$. By (\ref{eq-mixing}), we can choose some
$y_i\in U_{t_i}\cap f^{-N}U_{t_{i+1}}$. Then
$$d(f^{n_i}x_{N_i}^i,y)\leq diam U_{t_i}<\delta,\,\,d(f^N y_i, x_{1}^{i+1})\leq diam U_{t_{i+1}}<\delta.$$
By (\ref{eq-Yk-choice-2}), we know that
$$ d(f^{n_i}x_j^i,x_{j+1}^i)<\delta, j=1,2,\cdots,x_{N_i},i=1,2,\cdots,k.  $$
Then 
the  following  orbit segments form  a $\delta-$pseudo-orbit:
\begin{eqnarray*}
& & \cdots,\, f^{-t-1}x_1^1,\,f^{-t}x_1^1,\,f^{-t+1}x_1^1,\,\cdots,\, f^{-1} x_1^1,\,\\
& & x_1^1,\,f  x_1^1,\,\cdots,\, f^{n_1-1} x_1^1,\,\,\,\,\,\,\,\,x_2^1, \,\cdots,\, f^{n_1-1} x_2^1,\,\, \,\,\,\,\, \cdots,\,\,\,\,\,x_{N_1}^1,\,f  x_{N_1}^1,\,\cdots,\, f^{n_1-1} x_{N_1}^1,\,\\
& & y_1,\,fy_1,\,f^2y_1,\,\cdots,\,f^{N-1}y_1,\,\\
& & x_1^2,\,f  x_1^2,\,\cdots,\, f^{n_2-1} x_1^2,\,\, \,\,\,\,\,x_2^2, \,\cdots,\, f^{n_2-1} x_2^2,\,\,\,\,\,\,\,\,\,\,\, \cdots,\, \,\,\,x_{N_2}^2,\,f  x_{N_2}^2,\,\cdots,\, f^{n_2-1} x_{N_2}^2,\,\\
& & y_2,\,fy_2,\,f^2y_2,\,\cdots,\,f^{N-1}y_2,\,\\
& &\cdots,\\
& & x_1^i,\,f  x_1^i,\,\cdots,\, f^{n_i-1} x_1^i,\,\,\, \,\,\,\,\,\,x_2^i,\, \,\cdots,\, f^{n_i-1} x_2^i,\,\ ,\,\,\,\, \cdots,\,\,\,\,\,\,\,x_{N_i}^i, \cdots,\, f^{n_i-1} x_{N_i}^i,\,
\\& & y_i,\,fy_i,\,f^2y_i,\,\cdots,\,f^{N-1}y_i,\,\\
& &\cdots, \\
& & x_1^{k-1},\, fx_1^{k-1}, \,\cdots,\, f^{n_{k-1}-1} x_1^{k-1},\,\,\,\,\,\,\,\,
\cdots,\,\,\,\,\,\,\,\,\,x_{N_{k-1}}^{k-1},\, \,\cdots,\, f^{n_{k-1}-1} x_{N_{k-1}}^{k-1},\,\\
& & y_{k-1},\,fy_{k-1},\,f^2y_{k-1},\,\cdots,\,f^{N-1}y_{k-1},\,\\
& & x_1^k,\,f  x_1^k,\,\cdots,\, f^{n_k-1} x_1^k,\,\,\, \,\,\,\,\,\,x_2^k,\, \,\cdots,\, f^{n_k-1} x_2^k,\,\ ,\,\,\,\, \cdots,\,\,\,\,\,\,\,x_{N_k}^k, \cdots,\, f^{n_k-1} x_{N_k}^k,\,\\
  & & f^{n_k} x_{N_k}^k,\,\,f^{n_k+1} x_{N_k}^k,\,\cdots, f^{t } x_{N_k}^k,\,f^{t+1 } x_{N_k}^k,\cdots\\
\end{eqnarray*}
By exponential shadowing, there is a point $x$ whose orbit $(\tau,\lambda)-$shadows above pseudo-orbit and thus
$x\in B(\ul x_1, \ldots \ul x_k)$.

We define $F_k$ by
\[
F_k = \{\overline{B(\ul x_1, \ldots, \ul x_k)}: (\ul x_1, \ldots \ul x_k) \in S_1^{N_1} \times \ldots \times S_k^{N_k}\}.
\]
Note that $F_k$ is non-empty, compact and $F_{k+1} \subseteq F_k$. Define $F = \bigcap_{k=1}^{\infty} F_k$.  Remark that $F$ is compact and non-empty. 

\bl \label{2015Lemm-6}
For any $p \in F$, the sequence $\frac{1}{t_k} \log\|A(t_k,p)\|$ diverges, where $t_k = \sum_{i=0}^{k} N_i n_i+kN$.
\el
\bp
Choose $p \in F$ and let $p_k := f^{t_{k-1}}p$. Then there exists $(x_1^k, \ldots, x_{N_k}^k) \in S_k^{N_k}$ such that
\[
p_k \in \bigcap_{j=1}^{N_k} f^{-(j-1)n_k} \overline{\bigcap_{t=0}^{n_k-1}f^{-t}B ( f^t x^k_j, \tau e^{\lambda \min\{t,n_k-t\}})}.
\]

We follow the idea of \cite{Kal} to give following estimates.
Firstly let us consider $k$ with $k=1 \text{ mod } 2. $  Fix $\,j=1,2,\cdots,N_k$ and consider the orbit segments $$p_{k+(j-1)n_k},fp_{k+(j-1)n_k},\cdots,f^{n_k}p_{k+(j-1)n_k}$$ and $x_j^k,fx_j^k,\cdots,f^{n_k}x_j^k.$
Since the Lyapunov splitting and Lyapunov metric are continuous on the
compact set $\mathcal{R}^{\mu_k}_{\epsilon,l}$, the cones $K^\eta_0(x_{j+1}^k)\subseteq K_0=K_0(x_{j+1}^k)$ (centered at point $x_{j+1}^k$) and $K^\eta_{0}( f^{n_k}x_j^k)$ are close
if $x_{j+1}^k$ and $  f^{n_k}x_j^k$ are close enough. Therefore we can ensure that
$K^\eta_{0}( f^{n_k}x_j^k) \subset K_0(x_{j+1}^k)$ if $\delta$ small enough and thus
$A(p_{k+(j-1)n_k},n_k) (K_0(x_j^k)) \subset K_0(x_{j+1}^k)$. Using the norm estimate
in    Lemma \ref{Lem-simple-estimate-Lyapunov-2} (in this estimate $\chi=a$, being the largest Lyapunov exponent of $\mu_k=\mu_1$) and   (\ref{eq-different-norm-estimate}), (\ref{eq-u-u'}),  we obtain for any $u \in K^\eta_0(x_{j}^k)$
$$
\| A(p_{k+(j-1)n_k},n_k)\, u \| _{n_k} \ge \| (A(p_{k+(j-1)n_k},n_k)\, u)' \| _{n_k} \ge  e^{n_k(a - 2\epsilon)} \| u' \| _0
$$
$$\ge \frac1{\sqrt{2}} e^{n_k(a - 2\epsilon)} \| u \| _0 \ge \frac12 e^{n_k(a - 2\epsilon)} \| u \| _{n_k}.
$$
Since $A(p_{k+(j-1)n_k},n_k)\, u \in K_0(x_{j+1}^k)$ for any $u \in K^\eta_0(x_{j}^k)$, we can iteratively apply $A(p,n)$ and
use the inequality above to estimate the largest Lyapunov exponent at $p$
$$
  \|  A(p_k,N_kn_k) u \|_{n_k} \ge   \left( \left(\frac12 e^{n_k(a - 2\epsilon)}\right)^{N_k}\| u \| _{n_k} \right).
$$
Then
we have
\[
\frac{1}{{n_k N_k}}\|A(n_k N_k,p_k) \|   \geq \frac1{n_k}\log \frac12 +a-2\epsilon\geq a-3\epsilon.
\]
It follows that
\[
\liminf_{k \rightarrow \infty}\frac{1}{{n_k N_k}}\|A(n_k N_k,p_k) \|   \geq a-3\epsilon.
\]
We can use the fact that $\frac{n_k N_k}{t_k} \rightarrow 1$ to prove that for $k$ with $k= 1\text{ mod } 2, $
\begin{eqnarray}\label{eq-a-estimate}
\liminf_{k \rightarrow \infty}\frac{1}{{t_k}}\|A(t_k,p) \|   \geq\liminf_{k \rightarrow \infty}\frac{1}{{n_k N_k}}\|A(n_k N_k,p_k) \|   \geq a-3\epsilon.
\end{eqnarray}

Secondly let us consider $k$ with $k= 0 \text{ mod } 2. $  Fix $\,j=1,2,\cdots,N_k$, for the orbit segments $$p_{k+(j-1)n_k},fp_{k+(j-1)n_k},\cdots,f^{n_k}p_{k+(j-1)n_k}$$ and $x_j^k,fx_j^k,\cdots,f^{n_k}x_j^k,$ by the first estimate in (\ref{eq-estimate-simple-Lyapunov-new})
 of Lemma \ref{Lem-New-simple-estimate-Lyapunov-new} (in this estimate $\chi=b$, being  the largest Lyapunov exponent of $\mu_k=\mu_2$)   we have
$$\|A(n_k  ,\, p_{k+(j-1)n_k}) \|
\leq   l^2 e^{l}e^{n_k (b+\epsilon)} $$
Then
we have
\[
\|A(n_k N_k,p_k) \|\leq  \prod_{j=1}^{N_k} \|A(n_k  ,\, p_{k+(j-1)n_k}) \|
\leq   (l^2 e^{l})^{N_k}\cdot e^{n_k N_k(b+\epsilon)}
\]
and hence combing with $ l^2 e^{l}<e^{ n_2\epsilon}=e^{ n_k\epsilon}$  we have
\[
\frac{1}{{n_k N_k}}\|A(n_k N_k,p_k) \|   \leq \frac1{n_k}\log (l^2 e^{l}) + b+\epsilon\leq b+2\epsilon.
\]
It follows that
\[
\limsup_{k \rightarrow \infty}\frac{1}{{n_k N_k}}\|A(n_k N_k,p_k) \|   \leq b+2\epsilon.
\]
We can use the fact that $\frac{n_k N_k}{t_k} \rightarrow 1$ to prove that for $k$ with $k= 0 \text{ mod } 2, $
\begin{eqnarray}\label{eq-b-estimate}
\limsup_{k \rightarrow \infty}\frac{1}{{t_k}}\|A(t_k,p) \|   \leq\limsup_{k \rightarrow \infty}\frac{1}{{n_k N_k}}\|A(n_k N_k,p_k) \|   \leq b+2\epsilon.
\end{eqnarray}

Since the choice of $\epsilon$ satisfies $b+2\epsilon<a-3\epsilon,$ by  (\ref{eq-a-estimate}) and (\ref{eq-b-estimate}), the sequence $\frac{1}{{t_k}}\|A(t_k,p) \|$ diverges. \ep

\textbf{Step 3}.{Construction of a special sequence of measures $\omega_k$.} \label{ss}
We must first undertake an intermediate construction. For each $\ul x =(\ul x_1, \ldots, \ul x_k) \in S_1^{N_1} \times \ldots \times S_k^{N_k}$, we choose one point $z = z(\ul x)$ such that
\[
z \in B(\ul x_1, \ldots \ul x_k).
\]
Let $\mathcal  T_k$ be the set of all points constructed in this way. We show that points constructed in this way are distinct and thus $\#\mathcal  T_k = \#S_1^{N_1} \ldots \#S_k^{N_k}$. 
\bl \label{lb2}
Let $\ul x$ and $\ul y$ be distinct members of $S_1^{N_1} \times \ldots \times S_k^{N_k}$. Then $z_1 := z ( \ul x )$ and $z_2 := z( \ul y )$ are distinct points. Thus $\#\mathcal  T_k = \#S_1^{N_1} \ldots \#S_k^{N_k}$. 
\el
\begin{proof}
Since $\ul x \neq \ul y$, there exists $i,j$ so $x^i_j \neq y^i_j$. We have
\[
d_{n_i} (x^i_j, f^{h} z_1 ) < \tau<\varepsilon \mbox{ and }  d_{n_i} (y^i_j, f^{h} z_2 ) < \tau<\varepsilon,
\]
where $h= \sum_{l=0}^{i-1} N_l n_l +{i-1}N+ (j-1) n_i$. Since $d_{n_i} (x^i_j , y^i_j) >4 \varepsilon$,  we have $d_{n_i} (f^{h} z_1 , f^{h} z_2) >2 \varepsilon$.
\end{proof}

Now we start to define the measures on $F$ which yield the required estimates for the Entropy Distribution Principle.  For each $k$, define a Dirac measure centred on $\Tk$. More precisely, let
\[
\nu_k := \sum_{z \in \Tk}\delta_z
\]
We normalize $\nu_k$ to get a sequence of probability measures $\omega_k$, i.e.  $\omega_k := \frac{1}{\# \mathcal  T_k} \nu_k$.
\bl \label{2015Lemm-5}
Suppose $\omega$ is a limit measure of the sequence of probability measures $\omega_k$. 
Then $\omega (F) = 1$.
\el
\bp
For any fixed $l$ and   $\forall \,\,p \geq 0$, $\omega_{l+p} (F_{l}) = 1$ since $\omega_{l+p} (F_{l+p}) = 1$ and $F_{l+p} \subseteq F_{l}$.
Suppose $\omega = \lim_{k \rightarrow \infty} \omega_{l_k}$ for some $l_k \rightarrow \infty$, then $\omega(F_l) \geq \limsup_{k \rightarrow \infty} \omega_{l_k} (F_l) = 1$. It follows that $\omega (F) = \lim_{l \rightarrow \infty} \omega (F_{l}) = 1$.
\ep

\textbf{Step 4}. Estimate of $h_{top}^B ( MLI(A, f),f)$.

Let $\mathcal  B := B_n (q, \varepsilon)$ be an arbitrary ball which intersects $F$. Let $k$ be the unique number which satisfies $t_k \leq n < t_{k+1}$.

 We firstly consider $n$  with  $t_k \leq n < t_{k+1}-N.$ Let $j \in \{0, \ldots, N_{k+1} -1 \}$ be the unique number so
\[
t_k + n_{k+1}j \leq n < t_k + n_{k+1}(j+1).
\]
We assume that $j \geq 1$ and leave the details of the simpler case $j=0$ to the reader. The following lemma reflects the restriction on the number of points that can be in $\mathcal  B \cap \mathcal  T_{k+p}$.
\bl\label{2015Lemm-7}
For $p \geq 1$, $\omega_{k+p} ( \mathcal  B) \leq (\# \Tk )^{-1} (\#\mathcal  S_{k+1})^{-j}$
\el

\begin{proof} 
First we show that $\omega_{k+1} ( \mathcal  B) \leq (\# \Tk )^{-1} (\#\mathcal  S_{k+1})^{-j}$. We require an upper bound for the number of points in $\Tkk \cap \mathcal  B$.
If  $\omega_{k+1} (\mathcal  B) > 0$, then $\Tkk \cap \mathcal  B \neq \emptyset$. Let $z = z(\ul x, \ul x_{k+1}) \in \Tkk \cap \mathcal  B$ where $\ul x \in S_1^{N_1} \times \ldots \times S_k^{N_k}$ and $\ul x_{k+1} \in S_{k+1}^{N_{k+1}}$. Let
\[
\mathcal  A_{\ul x; x_1, \ldots, x_j} = \{ z(\ul x, y_1, \ldots, y_{N_{k+1}}) \in \Tkk : x_1 = y_1, \ldots, x_j = y_j \}.
\]
We suppose that $z^\prime = z(\ul y, \ul y_{k+1})\in \Tkk \cap \mathcal  B$ and show that $z^\prime \in \mathcal  A_{\underline{x}; x_1, \ldots, x_j}$. We have $d_n(z, z^\prime)<2 \varepsilon$ and we show that this implies $x_l = y_l$ for $l \in \{1, 2, \ldots, j\}$ (the proof that $\ul x = \ul y$ is similar). 
Suppose that $y_l \neq x_l$ and let $a_l = t_k + (l-1)(n_{k+1})$.  Using
\[
d_{n_{k+1}} (f^{a_l} z, x_l ) < \tau<{\varepsilon}  \mbox{ and }  d_{n_{k+1}} (f^{a_l} z^\prime, y_l ) <  \tau<{\varepsilon} .
\]
Recall that $d_{n_{k+1}} (x_l , y_l) >4 \varepsilon$.
Then we have
\beq
d_{n} (z, z^\prime) &\geq& d_{n_{k+1}} (f^{a_l}z, f^{a_l}z^\prime) \\
&\geq& d_{n_{k+1}} (x_l, y_l) - d_{n_{k+1}} (f^{a_l}z, x_l) - d_{n_{k+1}} (f^{a_l}z^\prime, y_l) >2\varepsilon,
\eeq
which is a contradiction.
Thus, we have
\beq
\nu_{k+1} (\mathcal  B) \leq \# \mathcal  A_{x; x_1, \ldots, x_j} = (\#S_{k+1})^{N_{k+1} - j},
\eeq
\[
\omega_{k+1} (\mathcal  B) \leq (\# \mathcal  T_{k+1})^{-1} (\#S_{k+1})^{N_{k+1} - j} = (\# \Tk )^{-1} (\#\mathcal  S_{k+1})^{-j}
\]
Now consider $\omega_{k+p} ( \mathcal  B)$. Arguing similarly to above, we have
\beq
\nu_{k+p} ( \mathcal  B)& \leq &\# \mathcal  A_{x; x_1, \ldots, x_j} (\# \mathcal  S_{k+2})^{N_{k+1}} \ldots (\# \mathcal  S_{k+p})^{N_{k+p}}
\eeq
The desired result follows from this inequality by dividing by $\# \mathcal  T_{k+p}$.
\end{proof}
By (\ref{eq-Yk-choice}), the first estimate in (\ref{eq-E-choice}) and the first estimate in (\ref{eq-E-choice-2}), we have 
\begin{eqnarray}\label{eq-TkSk}
\# \mathcal  T_k (\# \mathcal  S_{k+1})^j &\geq& \exp \{ (h^*-3 \gamma) (N_1 n_1 + N_2 n_2 + \ldots + N_k n_k +j n_{k+1}) \}\nonumber\\
&\geq& \exp \{ (h^*-3 \gamma)   (n-kN) \}= \exp \{ (h^*-3 \gamma)n \exp \{ ( 3 \gamma- h^*)kN\}\nonumber\\
&=&\exp \{ (h^*-4 \gamma) n \} \exp \{ n\gamma+( 3 \gamma- h^*)kN\} \nonumber\\
&\geq& \exp \{ (h^*-4 \gamma) n \} \exp \{ \sum_{i=1}^kn_i\gamma+( 3 \gamma- h^*)kN\} \nonumber\\
&\geq&   \exp \{ (h^*-4 \gamma) n \}.\nonumber
\end{eqnarray}
 \begin{eqnarray}\label{eq-tkbetweentk1-N000000} \omega_{k+p} ( \mathcal  B) \leq    (\# \mathcal T_{k-1} )^{-1} (\#\mathcal  S_{k})^{-j}\leq  \exp \{ -(h^*-4 \gamma) n\} .
  \end{eqnarray}

Now we
  consider $n$  with  $ t_{k+1}-N \leq n < t_{k+1}.$ By Lemma \ref{2015Lemm-7} and (\ref{eq-TkSk}), using $$ B_{t_{k+1}-N-1}(q,\varepsilon) \supseteq \mathcal  B,$$ we have that for $p \geq 1$,   $$\omega_{k+p} ( \mathcal  B) \leq \omega_{k+p} ( \mathcal  B_{t_{k+1}-N-1}(q,\varepsilon))\nonumber \\
     \leq    \exp \{ -(h^*-4 \gamma) ({t_{k+1}-N-1})\}.$$
 Combining this with the first estimate in (\ref{eq-E-choice}) and the first estimate in (\ref{eq-E-choice-2}), we have
    \begin{eqnarray}\label{eq-tk-Nbetweentk} \,\,\,\,\,\,\, \omega_{k+p} ( \mathcal  B) \leq    \exp \{ -(h^*-4 \gamma) ({t_{k+1}-N-1})\} \leq \exp \{ -(h^*-5\gamma) n \} .
  \end{eqnarray}


Combining (\ref{eq-tkbetweentk1-N000000}) and  (\ref{eq-tk-Nbetweentk})
we have for all $n$,
\begin{eqnarray}\label{eq-estimate-bowenboll-sup}
\limsup_{l \rightarrow \infty} \omega_{l} (B_n (q, \varepsilon)) \leq \exp \{-n(h^* - 5 \gamma)\}.
\end{eqnarray}

By (\ref{eq-estimate-bowenboll-sup}) and Entropy Distribution Principle (Proposition \ref{Prop2015-Pressure-distribution}), $h_{top}^B(F,\varepsilon,f)\geq h^* - 5\gamma.$
Since $\gamma$ and $\varepsilon$ were arbitrary and $F \subset   MLI(A, f)$, we have $$h_{top}^B ( MLI(A, f),f)  \geq h^*.$$

Remark that there is another way to get this estimate. More precisely,  from (\ref{eq-estimate-bowenboll-sup}) we know that $\omega(B_n (q, \varepsilon)) \leq \exp \{-n(h^* - 4 \gamma)\}.$ So for any $q\in F,$ $$\underline{h}_\mu(f,q)\geq h^* - 5 \gamma. $$ By (\ref{equation-bowenentropy})
$$h_{top}^B(F, f)=\sup \{\underline{h}_\mu(f): \mu \in \mathcal{M}(X), \; \mu(F)=1\}\geq h^* - 5 \gamma.$$

\textbf{Step 5}. Estimate of $h_{top}^P ( MLI(A, f),f)$.  By (\ref{eq-Yk-choice}), the first estimate in (\ref{eq-E-choice}) and the first estimate in (\ref{eq-E-choice-2}), we have
\beq
\# \mathcal  T_{k}    \geq  \# \mathcal  S_{k}    \geq \exp \{ (h_{\mu_k}(f)-3 \gamma)  N_k n_k   \}.
\eeq
Without loss of generality, assume that $H^*=h_{\mu_1}(f).$  Combining this with Lemma \ref{2015Lemm-7} ($n=t_k$) and $\lim_{k\rightarrow \infty} \frac{N_kn_k}{t_k}=1$, for large enough $k$ with $k=1\text{ mod }2,$ we have
\begin{eqnarray}\label{eq-estimate-bowenboll-inf}
\limsup_{l \rightarrow \infty} \omega_{l} (B_{t_k} (q, \varepsilon))\leq (\# \mathcal  T_{k})^{-1} \leq \exp \{-t_k(H^* - 4 \gamma)\}.
\end{eqnarray}
This implies that  for large enough $k$ with $k=1\text{ mod }2,$
$$\omega(B_{t_k} (q, \varepsilon)) \leq \liminf_{l \rightarrow \infty} \omega_{l} (B_{t_k} (q, \varepsilon))\leq \exp \{-{t_k}(H^* - 4 \gamma)\}.$$  So for every $q\in F,$ $$\overline{h}_\mu(f,q)\geq H^* - 4 \gamma. $$ Since $F$ is compact, by (\ref{equation-Packentropy})
$$h_{top}^P( F,f)=\sup \{\overline{h}_\mu(f): \mu \in \mathcal{M}(X), \; \mu(F)=1\}\geq H^* - 4 \gamma.$$
Since $\gamma$ and $\varepsilon$ were arbitrary and $F \subset   MLI(A, f)$, we have $h_{top}^P ( MLI(A, f),f)  \geq H^*$. Now we ends the proof of Proposition \ref{prop-entropy-twomeasure}. \qed

\begin{Rem}\label{Rem-anotherway} The  construction of $F$ in above proof is very cautious and mainly used to deal with $ h_{top}^B.$  If we only need to estimate   $h_{top}^P$ (and $h_{top}$) in  Step 5, we can choose another (little) simple constructed way: one just modify $N_k$ by $N_k\equiv1$ and  $n_k$ by a sequence with $ \lim_{k\rightarrow \infty}\frac{n_k}{n_{k+1}}=0.$ In other words, we only need to use every $\mathcal S_k$ just once. One can follow above proof straightforward to estimate   $h_{top}^P$ (and $h_{top}$) and here we omit the details. This is based on we only need a uniform estimate for  the measure of $B_n (q, \varepsilon)$ of a subsequence of $n=t_k$. However, we emphasize that the new way is not valid to deal with $ h_{top}^B,$ since we need a uniform estimate for the measure of   $B_n (q, \varepsilon)$ of arbitrarily large $n$.
\end{Rem}

Now we continue the proof of Theorem \ref{Thm-Entropy-Maximal-LyapunovIrregular-HolderCocycle}. By (1) of Proposition \ref{Prop-Cor2015-Contiuity-LyaExp} and Proposition \ref{prop-entropy-twomeasure},
$$h_{top}(f) =\sup_{\mu, \nu\in \mathcal{M}^e_{f}(X)} \{\max\{ h_\mu(f), h_\nu(f)\} |\,\,\chi_{max}(\mu,A)>  \chi_{max}(\nu,A)\}$$ $$\leq \min\{ h_{top}(MLI(A,f),f) , h_{top}^P(MLI(A,f),f)\}.$$
It follows that
$$h_{top}(f) = h_{top}(MLI(A,f),f)= h_{top}^P(MLI(A,f),f) .$$

It is not difficult to check that a topologically mixing system with shadowing satisfies specification property (c.f. \cite [Proposition 23.20, Definition 21.1]{DGS}).
Since $f$ is topologically mixing and has exponential shadowing (stronger than shadowing), then $f$ has specification and thus by      \cite[Theorem 2.1 ]{PS2005} (or see  \cite{PS}, also  see  \cite{EKW} for similar discussion), $f$ has entropy-density property. By (2) of Proposition \ref{Prop-Cor2015-Contiuity-LyaExp} and Proposition \ref{prop-entropy-twomeasure},  $$h_{top}^B(MLI(A,f),f)\geq\sup_{\mu, \nu\in \mathcal{M}^e_{f}(X)} \{\min\{ h_\mu(f), h_\nu(f)\} |\,\,\chi_{max}(\mu,A)>  \chi_{max}(\nu,A)\} $$ $$=\sup_{\mu \in \mathcal{M}^e_{f}(X)} \{ h_\mu(f) |\,\,\chi_{max}(\mu,A)>\inf_{\nu\in \mathcal{M}_f(X)} \chi_{max}(\nu,A)\}. $$
Now we complete the proof of Theorem \ref{Thm-Entropy-Maximal-LyapunovIrregular-HolderCocycle}. \qed
\bigskip

By Theorem \ref{Thm-Entropy-Maximal-LyapunovIrregular-HolderCocycle} and (2) of Proposition \ref{Prop-Cor2015-Contiuity-LyaExp}, we have a direct consequence.

\begin{Cor}\label{Cor-thm-Entropy-Maximal-LyapunovIrregular-HolderCocycle}
Let $(X,d)$ be a compact metric space and $f:X \mapsto X$ be a topologically mixing  homeomorphism with   exponential shadowing
 property and let $A:M\rightarrow GL(m,\mathbb{R})$ be a H$\ddot{\text{o}}$lder continuous matrix function. If   $\chi_{max}(A,\mu)$ is  lower semi-continuous with respect to $\mu \in \mathcal{M}_f(X),$ then \\
 (1) either all ergodic measures have same Maximal Lyapunov exponents w.r.t. $A$, or \\
 (2) the entropy of maximal
Lyapunov irregular set of $A$ satisfies that
 $$h_{top}^B({MLI(A,f)},f)=h_{top}({MLI(A,f)},f)=h_{top}^P({MLI(A,f)},f)=h_{top}(f).$$
 \end{Cor}

\subsection{Proof of Theorem \ref{Thm-LyapunovIrregular-HolderCocycle} and \ref{Thm-horseshoe-LyapunovIrregular-Cocycle}}

By Proposition \ref{Prop2015-LemShadow}, $f|_X$ has exponential shadowing.
 Then we can apply Theorem \ref{Thm-Entropy-Maximal-LyapunovIrregular-HolderCocycle} to   complete the proof of Theorem \ref{Thm-LyapunovIrregular-HolderCocycle}.
   Applying Theorem \ref{Thm-LyapunovIrregular-HolderCocycle} for cocycle $A(x,n)=D_xf^n$, one ends the proof of Theorem \ref{Thm-horseshoe-LyapunovIrregular-Cocycle}. \qed


\section*{Acknowlegements}   The research of X. Tian was  supported by National Natural Science Foundation of China  (grant no. 11301088) and  Specialized
  Research Fund for the Doctoral Program of Higher Education  (No.  20130071120026).


\end{document}